\documentclass{article}
\usepackage{amsmath,amssymb,url,latexsym,amsthm}
\usepackage[all]{xy}
\usepackage{color}

\newtheorem{defi}{Definition}[section]
\newtheorem{theo}{Theorem}[section]
\newtheorem{lemma}{Lemma}[section]
\newtheorem{prop}{Proposition}[section]
\newtheorem{coro}{Corollary}[section]

\def\into{ \rightarrowtail }
\def\onto{ \twoheadrightarrow }
\def\splito{ \rightleftarrows }

\newdir{ >}{{}*!/-8pt/@{>}}

\def\EE{ \mathbb{E} }
\def\CC{ \mathbb{C} }
\def\DD{ \mathbb{D} }
\def\VV{ \mathbb{V} }

\begin{document}

\author{Dominique Bourn}

\title{On congruence modular varieties\\ and Gumm categories}

\date{}

\maketitle

\section*{Introduction}

In \cite{Gu}, Gumm characterized the congruence modular varieties  among all the varieties of Universal Algebra by the validity of the \emph{Shifting Lemma}. Doing this, he translated the congruence modular formula:
$$(T\vee S) \wedge R=T\vee (S\wedge R), \rm{\; for \; any \; triple\;} (T,S,R)\; \rm{\; such \; that:\;} T\subset R$$
in geometric terms: given any triple of equivalence relations $(T,S,R)$ such that $R\cap S\subset T$, the following left hand side situation implies the right hand side one:
$$ 
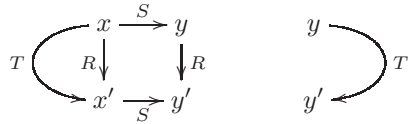
\xymatrix@=15pt{
	x \ar[r]^{S} \ar@(l,l)[d]_T \ar[d]_R & y \ar[d]^R &&  y \ar@(r,r)[d]^T\\
	x' \ar[r]_{S} & y' &&   y'
}
$$
One of the main interest of the Shifting lemma is that, being freed of any inductive condition (as the existence of suprema for instance), it allows us to recover in any category $\VV(\EE)$ of  internal $\VV$-algebras in $\EE$ all the projective properties of the variety $\VV$, when this variety is  congruence  modular.

Being non-inductive, this shifting property has a meaning in any finitely complete category $\EE$; its categorical characterization was given in \cite{BGGGum} where its first consequences were investigated, mainly concerning the centralization of equivalence relations, the internal groupoids and the internal categories. Later on, the categories satisfying the Shifting Lemma  were called \emph{Gumm categories} in \cite{BGGum} and \cite{BGum}. 

Now, in the same way as the Mal'tsev categories \cite{Bfib}, the Gumm categories were characterized in \cite{BGum} by some properties of the pointed fibers $Pt_Y\EE$ (=split epimorphisms above  the  object $Y$) of their fibrations of points $\P_{\EE}$, see Sections \ref{section1} and \ref{section2} below. The aim of this work is to start from this characterization to deepen the investigatation on the Gumm categories, to better understand the observations already made in \cite{BGGGum} and to produce new examples. Indeed, the characterization theorem will provide us with a golden thread inside the Gumm categories: an intrinsic notion of abelian split epimorphism, see Section \ref{absplit}.

A natural effect of these investigations will be to make it possible to measure the distance between the consequences of the Shifting Lemma in the varietal context and in the much more general categorical one: in Section \ref{axiomst} is made explicit a property (Axiom $\ast$) which is satisfied by the congruence modular varieties, but not in general by the Gumm categories, and whose main interest, when it holds, is to produce the universal abelian split epimorphism, see Theorem \ref{main2}. Now, beyond that, these  investigations draw  attention on a very specific phenomenon which we are going to describe and  emphasize below,  before entering into the details of the organization of this article.

\medskip

\noindent\textbf{Algebraic chrystallography}

Let us define as \emph{chrystallographic for a given algebraic structure} any categorical setting in which, on any object $X$ of this setting, there is at most one internal algebraic structure of this  kind.

Such kinds of situations are well known from a long time; for instance, it is clear that in a pointed J\'onnson-Tarski variety, on any algebra $X$ there is at most one internal commutative  monoid structure; the same property holds for the commutative  and associative (=autonomous) Mal'tsev operations in the  Mal'tsev varieties. The first property still holds for the unital categories \cite{Bint} and the second one for the Mal'tsev categories \cite{CLP} and \cite{CPP}. We can even widen the range of examples by saying that the setting $RGh\EE$ of internal reflexive graphs in a Mal'tsev category $\EE$ is chrystallographic for the  notion of internal groupoid in $\EE$.

But, in a way, the particular varietal origin of the first examples lessened the surprise effect: this phenomenon arised naturally because, and when, some term of the variety $\VV$ became a homomorphism.

Now, along with some aspects of our investigations, emerged a situation where  it is impossible to skip the surprise effect. Consider any pointed variety $\VV$ (i.e. with only one constant $0$) endowed with one binary  term $\circ$ and one unary term $\upsilon$ satisfying:
\begin{align*}
a\circ 0=0\circ a;\;\;\;\;\;\; ; \;\;\ & \;\;\;\; \upsilon(a\circ 0)=a 
\end{align*}
This seems to be a very poor  varietal context. However, on any algebra $X$ in such a variety, there is at most one internal subtraction (which is then the subtraction associated  with an abelian group structure), see Proposition \ref{abob}. And, in this context, this subtraction is far from being (at least directly) determined by the terms of the variety in question.

\medskip

Here is the organization of this article. After a first section devoted to conventions and notations, Section 2 recalls the Characterization Theorem and investigates the three ways opened, within the pointed categories, by this characterization. The  most important one seems to be the setting of \emph{congruence hyperextensible categories} whose one example is given by the varieties with the terms described above, and whose main structural property is to be crystallographic for the subtractions and thus to determine an intrinsic notion of abelian object. Section 3 describes, with the notion of \emph{diagonal punctuation}, i.e. a way to produces universal abelian  objects in the pointed congruence modular varieties. Section 4 investigates the  structural outcomes of the Characterization Theorem in following the golden thread provided by the notion of abelian split epimorphisms which are nothing but the intrinsic abelian objects of the congruence hyperextensible  fibers $Pt_Y\EE$. Sections \ref{cat1} and \ref{cat2} re-investigate and enlarge the results of \cite{BGGGum} on the internal groupoids and internal categories, producing important structural precisions and new classes of examples of Gumm categories. Section 5 is devoted to Axiom $\ast$ which allows to recover  in an abstract way the results of Section 3 concerning the relationship between  the diagonal punctuation and the construction of the universal abelian split epimorphism. Finally Section 6 collects all the needed informations about the  relationship between congruence permutation and  centralization of pairs of congruences.

\section{Conventions, notations}\label{section1}

In this article, any category $\EE$ will be supposed finitely complete.
Given any map $f:X\to Y$, we use the following simplicial notations for its kernel equivalence relation $R[f]$:
$$\xymatrix@=15pt
{
	R_2[f]  \ar@<2ex>[rr]^{d_{2}^f} \ar[rr]|{d_{1}^f} \ar@<-2ex>[rr]_{d_{0}^f}\ar[dd]_{R_2(g)} && R[f] \ar@<2ex>[rr]^{d_{1}^f} \ar@<-2ex>[rr]_{d_{0}^f} \ar[dd]_{R(g)} &&
	X   \ar[rr]^{f} \ar[ll]|{s_0^f} \ar[dd]_{g} && Y \ar[dd]^{h}\\
	&&&&\\
		R_2[f']  \ar@<2ex>[rr]^{d_{2}^{f'}} \ar[rr]|{d_{1}^{f'}} \ar@<-2ex>[rr]_{d_{0}^{f'}} && R[f'] \ar@<2ex>[rr]^{d_{1}^{f'}} \ar@<-2ex>[rr]_{d_{0}^{f'}} &&
	X'   \ar[rr]_{f'} \ar[ll]|{s_0^{f'}}  && Y' 
}
$$
and the notation $R(g)$ for the factorization induced by a right hand side commutative square.
A \emph{regular} epimorphism $f$ is a map $f$ which is the coequalizer of the two  legs of it kernel equivalence relation, namely its \emph{quotient map}.
More generally we use the simplicial notations for any internal category as well:
$$\xymatrix@=22pt
{
	X_{\bullet}:\;\;\;\;\;\;\;\;\;\;	X_1\times_0X_1  \ar@<2ex>[rr]^>>>>>>>>>{d_{2}^{X_{\bullet}}} \ar[rr]|>>>>>>>>>{d_{1}^{X_{\bullet}}} \ar@<-2ex>[rr]_>>>>>>>>>{d_{0}^{X_{\bullet}}} && X_1 \ar@<2ex>[rr]^{d_{1}^{X_{\bullet}}} \ar@<-2ex>[rr]_{d_{0}^{X_{\bullet}}} &&
	X_0   \ar[ll]|{s_0^{X_{\bullet}}}  
}
$$
In particular, any preorder $T$ on an object $Y$ will be denoted:
$$\xymatrix{  T\times_YT \ar@<-10pt>[rr]_{d_0^T}\ar[rr]|{d_1^T}  \ar@<10pt>[rr]^{d_2^T} 	&& T \ar@<-10pt>[rr]_{d_0^T} \ar@<10pt>[rr]^{d_1^T} && Y \ar[ll]|{s_0^T} 
}
$$
in other words the maps $s_0^T: Y\into X$ and $d_1^T:T\times_YT \to T$ will respectively denote the reflexivity and  transitivity morphisms.

Given any category $\EE$, the slice category $\EE/Y$ is the category whose objects are the maps with codomain $Y$ and  morphims the commutative triangles above $Y$, while the coslice category $Y/\EE$ is the category whose objects are the maps with domain $Y$ and morphisms are the commutative triangles under $Y$. As we shall see, specially important are the properties of $\EE$ which are stable under slicing and coslicing. 

A category $\EE$ is \emph{regular}, see \cite{Barr}, when the first two following properties hold and \emph{exact} when the three properties hold:\\
1) the regular epimorphisms are stable under pullback;\\
2) any kernel equivalence relation $R[f]$ admits a quotient map;\\
3) any equivalence relation $R$ is effective, i.e. of the form $R[f]$ for some map $f$.

Any variety $\VV$ is exact and a regular  epimorphism is nothing but a surjective homomorphism. Consider now any commutative quadrangle of regular epimorphisms in a regular category $\EE$:
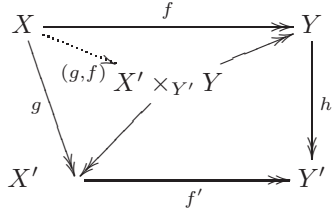
$$\xymatrix@=8pt
{
	X   \ar@{->>}[rrrr]^{f} \ar@{->>}[dddr]_{g} \ar@{.>}[drr]_{(g,f)}&&&& Y \ar@{->>}[ddd]^{h}\\
	&& X'\times_{Y'}Y \ar@{->>}[rru]_{} \ar@{->>}[ddl]_{}\\
	&&&&\\
X'  & \ar@{->>}[rrr]_{f'}   &&& Y' 
}
$$
When the factorization $(g,f)$ towards the pullback of $h$ along $f'$ is a regular epimorphism as well, the square in question is necessarily a pushout; we shall need the following:
\begin{defi}
Let $\EE$ be a regular category. A commutative square of regular epimorphims as above with a regular epic factorization $(g,f)$ is called a regular pushout.	
\end{defi}
The main interest of this notion is due to the fact that any regular pushout insures that the factorization $R(g)$ and all the factorizations $R_n(g)$ in our first diagram are regular epimorphisms as well.

\subsection{The fibration of points}

We denote by $Pt\mathbb{E}$ the category whose objects are the split epimorphisms in $\mathbb{E}$ with a given splitting and morphisms the commutative squares between these data; by $\P_{\mathbb E} :Pt\mathbb{E}\rightarrow \mathbb{E}$ we denote the functor associating its codomain with any split epimorphism. It is left exact and a fibration whose cartesian maps are the pullbacks of split epimorphisms; it is called the \emph{fibration of points} \cite{B0}. It has an important classification power in Algebra, see \cite{BB}. Let us recall two of them we shall need here \cite{Bfib}:
\begin{theo}\label{class1}
	Given any category $\EE$:\\
	1) it is a Mal'tsev category if and only if any fiber $Pt_Y\EE$ is unital;\\
	2) it is a naturally Mal'tsev category if and only if any fiber $Pt_Y\EE$ is additive.
\end{theo}
A category is a Mal'tsev one if and only if any reflexive relation is an equivalence relation, see \cite{CLP} and  \cite{CPP}. It is a naturally Mal'tsev category when any object $X$ is endowed with a natural Mal'tsev operation, see \cite{Jo}. A unital category fullfils the categorical characterization of the J\'onsson-Tarsky varieties, namely it is a pointed category $\EE$ in which the only punctual relation between $X$ and $Z$, see next Section, is given by the product $X\times Z$, see \cite{Bfib}. The more usual example is produced by the category $Mon$ of monoids.  
We shall recall, in Section \ref{kGumm} below, how this  same fibration $\P_{\mathbb E}$ characterizes the congruence modular varieties and the  Gumm categories as well \cite{BGum}.
The fibre above an object $Y\in \EE$, denoted by $Pt_Y\mathbb E$, is nothing but the coslice of a slice category, namely: $1_Y/(\EE/Y)$.

\subsection{Punctual spans in pointed categories}

A category $\EE$ is pointed when the terminal object $1$ is initial as well.  Any fiber $Pt_Y\EE$ is pointed. In a pointed category $\EE$, a span $(f,g)$:
$$\xymatrix@=4pt{
	& W  \ar@<-1ex>[ddd]_{f} \ar@<1ex>[rrrr]^{g}  &&&&  Z  \ar@<-1ex>[ddd]_{\tau_Z} \ar@{.>}[llll]^{t} \\
	&&&&& && \\
	&&& &&&&&&\\
	& X  \ar[rrrr]^{\tau_X}\ar@{.>}[uuu]_{s} &&&& 1 \ar@<1ex>[llll]^{0_X} \ar[uuu]_{0_Z}
}
$$
is said to be \emph{right punctual} (resp. \emph{left punctual}) when there is a section $t$ of $g$ (resp. $s$ of $f$) such that $f.t=0$ (resp. $g.s=0$). It is said to be punctual when it is both right and  left punctual. A relation is a span such that $(f,g): W\to X\times Z$ is a monomorphism.

Later on, we shall be highly interested in the punctual spans in the pointed fibers $Pt_Y\EE$, namely in the following commutative squares of split epimorphisms in $\EE$:
$$\xymatrix@=4pt{
	& W  \ar@<-1ex>[ddd]_{f} \ar@<1ex>[rrrr]^{g}  &&&&  Z  \ar@<-1ex>[ddd]_{\phi_Z} \ar[llll]^{t} \\
	&&&&& && \\
	&&& &&&&&&\\
	& X  \ar[rrrr]^{\phi_X}\ar[uuu]_{s} &&&& Y \ar@<1ex>[llll]^{\sigma_X} \ar[uuu]_{\sigma_Z}
}
$$
\begin{prop}
	Given any punctual span in the fiber $Pt_Y\EE$,  the idempotent morphisms $sf$ and $tg$ on $W$ do commute in $\EE$.
\end{prop}
\proof
The proof is straightforward: $sftg=s\sigma_X\phi_zg=t\sigma_Z\phi_Xf=tgsf$.
\endproof

So, the whole diagram originates from a pair commutative idempotents on the object $W$.

\subsection{Subtractions, opsubtractions and groups}\label{neg}

\begin{defi}
A binary operation $s$ on a set $X$ is said to be:\\
-1)  annihilating when there is an element $0$ such that $s(x,x)=0$ for all $x\in X$;\\
-2) a subtraction when $s(x,0)=x$ and $s(x,x)=0$\\
-3) an opsubtraction when $s(0,x)=x$, and $s(x,x)=0$.
\end{defi}
Obviously, any group structure $G=(X,\circ,0)$ provides a subtraction  with the \emph{difference mapping} $d_G(u,v)=u\circ v^{-1}$ and an opsubtraction with the \emph{opdifference mapping} $\partial_G(u,v)=u^{-1}\circ v$. Note that we get:
$$(*):\;\;\;\;\; \partial_G(u,v)\circ \partial_G(v,w)=\partial_G(u,w)$$
and thus $\partial_G$ produces an internal functor $\nabla_G\to G$:
$$\xymatrix@=4pt{
	 G\times G  \ar@<-2ex>[ddd]_{p_1^G} \ar@<2ex>[ddd]^{p_0^G}\ar[rrrr]^{\partial_G} &&&& G \ar@<-1ex>[ddd]_{\tau_G} \\
	&&&&& && &&\\
	&&& &&&&&&&&\\
	G  \ar[rrrr]_{\tau_G}\ar[uuu]|{s_0^G} &&&&  1 \ar[uuu]_{0}
}
$$
With the difference mapping of a group such that $x\circ x=0$, we have an example of a subtraction which is an opsubtraction as well. The following observation is well known probably from Tarski \cite{Ta}:

\begin{prop}\label{group}
	Let $(X,d,0)$  be any subtraction. The following condition are equivalent:\\
	1) $(X,d,0)$ is the subtraction associated with the difference mapping of a group;\\
	2) we have: $d(d(x,z),d(y,z))=d(x,y)$.
\end{prop}
\proof
Starting with a group $G=(X,\circ,0)$, the point 2) is straightforward.
Conversely, set $x\circ y=d(x,d(0,y))$. Then we get:
\begin{align*}
d(x,y)\circ y=&d(d(x,y),d(0,y))=d(x,0)=x\\
d(0,d(0,x))=&d(d(x,x),d(0,x))=d(x,0)=x\\
d(x\circ y,y)=&d(d(x,d(0,y)),y)=d(d(x,d(0,y)),d(0,d(0,y)))=d(x,0)=x
\end{align*}
In this way, $-\circ y$ becomes the inverse of $d(-,y)$; and we get $t\circ y =z \iff d(z,y)=t$. So:
\begin{align}
d(d(x,z),d(y,z))=d(x,y) \iff & d(x,y)\circ d(y,z)=d(x,z)\\
d(x,y)\circ d(y,z)=d(x,z) \iff &  t\circ d(y,z)=d(t\circ y,z)\\
x\circ d(y,z)=d(x\circ y,z) \iff &  x\circ (t \circ z)=(x\circ t) \circ z
\end{align}
(1) is obtained by composition of the left hand side equality with $-\circ d(y,z)$; (2) is obtained by setting $t=d(x,y)$ in the left hand side equality; (3) is obtained by setting $y=t\circ z$ in the left hand side equality, which is equivalent to $d(x\circ (t\circ z),z)=d((x\circ t)\circ z,z)$. 
The last assertion of the proposition is straighforward.
\endproof
We have the following precision \cite{BZ}:
\begin{prop}\label{abgroup}
	Let $(X,d,0)$ be the subtraction associated with a group structure. The following conditions are equivalent:\\
	1) the group is abelian;\\
	2) we have $d(x,d(x,y))=y$;\\
	3) $d$ satisfies the \emph{interchange law}:
	$ d(u,v)=d(u',v') \iff d(u,u')=d(v,v')$.
\end{prop}

\subsection{Ternary operation and fibration of points}

In any fibre $Pt_X\EE$ there is, with the split epimorphism $(p_1^X,s_0^X): X\times X\splito X$, an object which has a remarkable property with respect to the fibration $\P_{\EE}$. It is described by the following diagram where the right hand side part is a pullback: 
$$\xymatrix@=4pt{
	{X\;} \ar@{>->}[rrrr]^{(1_X,f)} \ar@<-1ex>[ddd]_{f}  &&&&  X\times Y  \ar@<-1ex>[ddd]_{p_1^Y} \ar@{ >->}[rrrr]^{X\times s} &&&& X\times X\ar@<-1ex>[ddd]_{p_1^X} \\
	&&&&& && &&\\
	&&& &&&&&&&&\\
	Y  \ar@{=}[rrrr]\ar[uuu]_{s} &&&& Y  \ar[uuu]_{(s,1_Y)}  \ar@{ >->}[rrrr]_{s} &&&& X \ar[uuu]_{s_0^X}
}
$$
and  means, among other things, that any split epimorphims with domain $X$ is ''embedded" in it, see \cite{Bmon}.
\begin{lemma}\label{biter}
	Given any category $\EE$, a binary operation on $(p_1^X,s_0^X)$ inside the fiber $Pt_X\EE$ is nothing but a ternary $p$ operation on $X$ in $\EE$ such that $p(x,x,x)=x$. This binary operation is a subraction in $Pt_X\EE$ if and only if it satisfies the Mal'tsev identities, $p(x,y,y)=x=p(y,y,x)$. This
	subtraction is the one associated with the difference mapping of a group in $Pt_X\EE$ if and only if the Mal'tsev law is associative. In this case the group is abelian if and only if the Mal'tsev law is autonomous. 
\end{lemma}
\proof
Straightforward is the first assertion: such a binary operation is necessarily of the form $d((x,z),(y,z))=(p(x,y,z),z)$; being pointed by $s_0^X$ in the fibre, we get $d((z,z),(z,z))=(z,z)$, namely $p(z,z,z)=z$.

We get  $d((x,z),(x,z)=(z,z)$ if and  only if we have $p(x,x,z)=z$. We get $d((x,z),(z,z))=(x,z)$ if and only if we have $p(x,z,z)=x$. See \cite{B3} for the  two last points.
\endproof
With the above ''embedding", we can now associate  a structural theorem for the fibers $Pt_Y\VV$ of any variety $\VV$:
\begin{theo}\label{teo}
	Let $\VV$ be a variety endowed with a ternary term $p$ satisfying $p(x,x,x)=x$. Then any object $(f,s): X\splito Y$ in the fiber $Pt_Y\VV$ is endowed with a pair $(\circ,d)$ of binary operations defined, for any pair $(u,v)\in R[f]$, by $u\circ v=p(u,sf(u),v)$ and $d(u,v)=p(u,v,sf(u))$.
\end{theo}
\proof
1) The base-change along $s$ in the above diagram determines, via the previous lemma, two binary operations on the split epimorphism $(p_1^Y,(s,1_Y))$ which are defined by:
$$(u,y)\circ (v,y)=(p(u,s(y),v),y) \;\; {\rm and} \;\; d((u,y),(v,y))=(p(u,v,s(y)),y)$$ 
2) Now, when $f$ is a $\VV$-homomorphism, these two operations are stable under the above inclusion $(1_X,f)$ in $Pt_Y\VV$; for we have for any pair $(u,v) \in R[f]$:
$$(u,f(u))\circ (v,f(v))=(p(u,sf(u)),v),f(u))=(1_X,f)(p(u,sf(u),v))$$
$$d((u,f(u)),(v,f(v)))=(p(u,v,sf(u)),f(u))=(1_X,f)(p(u,v,sf(u)))$$
\endproof

\section{Pointed variations on the Shifting Lemma}

\subsection{The characterization theorem}\label{section2}

As recalled in the introduction a Gumm category is a category satisfying the Shifting Lemma. Clearly, and importantly, this property is stable under slicing and coslicing, and consequently under the passage to the fibers $Pt_Y\EE$. We have also the following:
\begin{lemma}[\cite{BGGum}]
Let $F:\mathbb C\to \mathbb D$ be left exact and conservative functor (i.e. reflecting the isomorphisms). If $\mathbb D$ is a Gumm category, so is $\mathbb C$.
\end{lemma}
In particular, when $\DD$ is a Gumm category and $(T,\lambda,\mu)$  a left exact monad on $\DD$, then, the forgetful functor $U^T: Alg T \to \DD$ being left exact and conservative, the category $Alg T$ of $T$-algebras is a Gumm category as well.

Thanks to the Yoneda embedding any result is \cite{Gu} which is proved from the Shifting Lemma by using only projective properties remains valid in any Gumm category. It is the case, in particular, for the Cube Lemma (Proposition 2.4 in \cite{Gu}) we shall need later on:
\begin{prop}[Cube Lemma]\label{cubl}
Let $\EE$ be a  Gumm category. Given any triple of equivalence relations $(T,S,R)$ on an object $X$ such that $R\cap S\subset T$, the plain arrows imply the dotted one:
$$
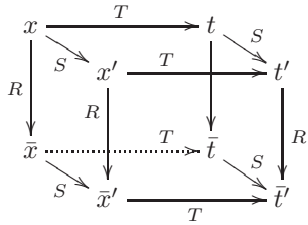
\xymatrix@=5pt{
	x \ar[rrd]_{S}  \ar[ddd]_{R}   \ar[rrrrr]^{T} &&&&& t \ar[rrd]^{S}  \ar[ddd]^{}\\
	&& x' \ar[ddd]_<<<<{R}  \ar[rrrrr]^<<<<<<{T}  &&&&& t' \ar[ddd]^{R}  && \\
	&&& &&&&&&\\ 
	\bar x \ar[rrd]_{S}  \ar@{.>}[rrrrr]^>>>>>{T}   &&&&& \bar t  \ar[rrd]^S  \\
	&&  \bar x'  \ar[rrrrr]_{T}    &&&&& \bar t' 
}
$$ 
\end{prop}

As recalled above, the fibration of points $\P_{\EE}$  allows to characterize the Gumm categories:
\begin{theo}[\cite{BGum}]
Given a category $\EE$, the following conditions are equivalent:\\
1) $\EE$ is Gumm category;\\
2) any fiber $Pt_Y\EE$ is congruence hyperextensive;\\
3) any fiber $Pt_Y\EE$ is congruence hypoextensive;\\
2) any fiber $Pt_Y\EE$ is punctually congruence modular.
\end{theo}

Any of these notions are recalled and investigated in the next sections. 

\subsection{Congruence hyperextensible categories}

\begin{defi}
	A pointed category $\EE$ is said to be congruence hyperextensible when given any punctual span in $\EE$ and any equivalence relation $T$ on $W$ such that $R[f]\cap R[g]\subset T$, we get $R[f]\cap g^{-1}(t^{-1}(T))\subset T$.  
\end{defi}
In other words, when the following situation holds in $W$:
$$ \xymatrix@=15pt{
	x \ar[r]^{R[g]} \ar@(l,l)@{.>}[d]_{T} \ar[d]_{R[f]} & tg(x) \ar[d]^{R[f]}  \ar@(r,r)[d]^{T}\\
	x' \ar[r]_{R[g]} & tg(x') 
}
$$
Any J\'onsson-Tarsky variety $\VV$ and  more generally any regular unital category $\EE$ is congruence hyperextensible \cite{BGum}. Here is an example of variety  which is clearly weaker than a J\'onsson-Tarski variety which is  obtained when we have $\upsilon=Id$:

\begin{prop}\label{exhyper}
	Let $\VV$ be a pointed variety with a binary term $\circ$ and one unary term $\upsilon$ such that:
	\begin{align*}
	a\circ 0=0\circ a;\;\;\;\;\;\; ; \;\;\ & \;\;\;\; \upsilon(a\circ 0)=a 
	\end{align*}
	Then the variety $\VV$ is congruence hyperextensible.
\end{prop}
\proof
Since any variety is an exact category, we can reduce the proof to punctual relations. So, consider any punctual relation between $X$ and $Z$ in $\VV$:
$$\xymatrix@=4pt{
	& W  \ar@<-1ex>[ddd]_{\pi_X} \ar@<1ex>[rrrr]^{\pi_Z}  &&&&  Z  \ar@<-1ex>[ddd]_{\tau_Z} \ar[llll]^{0_Z^W} \\
	&&&&& && \\
	&&& &&&&&&\\
	& X  \ar[rrrr]^{\tau_X}\ar[uuu]_{0_X^W} &&&& 1 \ar@<1ex>[llll]^{0_X} \ar[uuu]_{0_Z}
}
$$
and an equivalence relation $T$ on $W$ such that $(0Wb)T(0Wc)$. Then we get: 
\begin{align*}
(aW0)\circ(0Wb)&=(a\circ 0)W(0\circ b)=(a\circ 0)W(b\circ 0)& =(aWb)\circ(0W0)\;\;\;\;\;\; \\
 \;\;{\rm whence:}\;\; &\;\;\;((aWb)\circ(0W0))T((aWc)\circ(0W0)) &
\end{align*}
Then, applying $\upsilon$ to each side, we get: $(aWb)T(aWc)$.
\endproof

We shall need a name for the situation introduced by the previous lemma:
\begin{defi}\label{VD}
	Let us call \emph{Hyperex} and denote by $\mathrm{Hex}$ the pointed variety defined by the only two previous terms and equations.
\end{defi} 
Thanks to the Yoneda embedding, given any category $\EE$, the category $\mathrm{Hex}(\EE)$ of internal $\mathrm{Hex}$-algebras in $\EE$ is congruence hyperextensible. 
\begin{lemma}\label{cokernel}
	Let $\EE$ be any congruence hyperextensible category. Consider any punctual relation:
	$$\xymatrix@=4pt{
		& W  \ar@<-1ex>[ddd]_{\pi_X} \ar@<1ex>[rrrr]^{\pi_Z}  &&&&  Z  \ar@<-1ex>[ddd]_{\tau_Z} \ar[llll]^{0_Z^W} \\
		&&&&& && \\
		&&& &&&&&&\\
		& X  \ar[rrrr]^{\tau_X}\ar[uuu]_{0_X^W} &&&& 1 \ar@<1ex>[llll]^{0_X} \ar[uuu]_{0_Z}
	}
	$$
	Then $\pi_Z$ is the cokernel of $0_X^W$. Moreover, it is the unique  map $h:W\to Z$ in $\EE$ such that: $(1): h.0_Z^W=1_Z$ and $(2): h.0_X^W=0$.
\end{lemma}
\proof
By the Yoneda embedding it is enough to check our assertion in the category $Set$ of sets. Let $f:W\to V$ such that $f.0_X^W=0$, namely $f(xW0)=0=f(0W0)$. So, we are in the following situation:
$$\xymatrix@=7pt{
	xWz \ar[dd]_{R[\pi_X]} \ar[rr]_{R[\pi_Z]}  \ar@(u,u)@{.>}[rr]_{R[f]} && 0Wz\ar[dd]^{R[\pi_X]} \\
	&&&\\
	xW0 \ar@(d,d)[rr]^{R[f]} \ar[rr]^{R[\pi_Z]} && 0W0
}
$$
whence, by the congruence hyperextensibility, the bended upper arrow since we have: $R[\pi_X]\cap R[\pi_Z]=\Delta_W \subset R[f]$. So, we get $f(xWz)=f(0Wz)=f(x'Wz)$, and $f$ annihilates the kernel  equivalence relation $R[\pi_Z]$. Since $\pi_Z$ is a split epimorphism, we get the desired factorization: $f.0_Z^W: Z\to V$.

Given any map $h: W\to Z$ satisfying $(2)$, we get a (unique) map $\phi:Z\to Z$ such that $\phi.\pi_Z=h$. From $(1)$, by composition with $0_Z^W$, we get $\phi=1_Z$ and $h=\pi_Z$.
\endproof
The previous lemma holds, in particular, when the square in question is a pullback, and allows to make explicit a very important property which is shared by the (conceptually stronger) unital categories:
\begin{theo}\label{ffll}
	Let $\EE$ be a congruence hyperextensible category. Given any object $Y$, the base-change functor $\EE \to Pt_Y\EE$ along the terminal map is fully faithfull.
\end{theo}
\proof
It is clearly faithful since $p_1^Z$ is a (split) epimorphism. Now consider the following diagram:
	$$\xymatrix@=4pt{
	X\times Z \ar[rdd]^>{(p_0^X,\phi)} \ar@<-1ex>[dddd]_{p_0^X} \ar@<1ex>[rrrr]^>>>>>>>{p_1^Z}  &&&&  Z \ar@{.>}[rdd]^{\psi} \ar@<-1ex>[dddd]_{} \ar[llll]^<<<<<<<{\iota_1^Z} \\
	&&&&& &&& \\
	& X\times Z'  \ar@<-1ex>[dddd]_{} \ar@<1ex>[rrrr]^{}  &&&&  Z'  \ar@<-1ex>[dddd]_{} \ar[llll]^{} \\
	&&&&& &&& \\
	X  \ar[rrrr]^{}\ar[uuuu]_{\iota_0^X} \ar@{=}[rdd] &&&& 1 \ar@{=}[rdd] \ar@<1ex>[llll]^{} \ar[uuuu]_{}\\
	&&& &&&&&&\\
	& X  \ar[rrrr]^{\tau_Y}\ar[uuuu]_{} &&&& 1 \ar@<1ex>[llll]^{0_Y} \ar[uuuu]_{0_{Z'}}
}
$$
Given any map $(p_0^X,\phi):X\times Z \to X\times Z'$ making commute the vertical left hand side quadrangle, the previous lemma determines a unique factorization $\psi: Z\to Z'$  necessarily making the upper horizontal quadrangle a pullback, which means: $(p_0^X,\phi)=X\times \psi$.
\endproof
\begin{prop}\label{abob}
	Let $\EE$ be a congruence hyperextensive category. Any subtraction $d$ on an object $X$ is the difference mapping associated with an internal group structure on $X$ which is necessarily abelian. On any object $X$, there is at most one subtraction $d$; any morphism $f:X\to Y$ between two objects equipped with a subtraction is a subtraction homomorphism. This subtraction is an opsubtraction if and only if the associated abelian group $(X,\circ,0)$ is such that $x\circ x=0$.
\end{prop}
\proof
Let $d$ be any subtraction on an object $X$. Consider the following punctual relation:
	$$\xymatrix@=4pt{
	& X\times X \times X  \ar@<-1ex>[ddd]_{p_0^X} \ar@<1ex>[rrrrr]^>>>>>>>>{p_1^X.p_2^X}  &&&&&  X  \ar@<-1ex>[ddd]_{\tau_X} \ar[lllll]^<<<<<<<{0_{X\times X}\times X} \\
	&&&&& && \\
	&&& &&&&&&\\
	& X\times X  \ar[rrrrr]^>>>>>>>>>>{\tau_{X\times X}}\ar[uuu]_{X\times X \times 0_X} &&&&& 0 \ar[uuu]_{0_X} \ar@<1ex>[lllll]^<<<<<<<<<<{0_{X\times X}}
}
$$
Setting $\beta(x,y,z)=d(d(x,z),d(y,z))$, we get the following situation associated with the previous punctual relation:
$$\xymatrix@=7pt{
	(x,y,z) \ar[dd]^{R[p_0]} \ar[rr]^{R[p_1.p_2]}  \ar@(l,l)@{.>}[dd]_{R[\beta]} && (0,0,z)\ar[dd]_{R[p_0]} \ar@(r,r)[dd]^{R[\beta]}\\
	&&&\\
	(x,y,0)  \ar[rr]_{R[p_1.p_2]} && (0,0,0)
}
$$
since we have the right hand side bended arrow by $d(t,t)=0$. When $\EE$ is congruence  hyperextensive, we get: $\beta(x,y,z)=\beta(x,y,0)$, namely $d(d(x,z),d(y,z))=d(d(x,0),d(y,0))=d(x,y)$. So we get the first point with the group law given by $x\circ y=d(x,d(0,y))$. The inverse of $x$ is $d(0,x)$; so we get $(*): d(0,d(0,x))=x$.

We know that the group in question is commutative if and  only if we have $d(x,d(x,y))=y$. Then define $\psi: X\times X \to X$ by $\psi(x,y)=d(x,d(x,y))$, so that, by $(*)$, we have $\psi(0,y)=y$; we also get: $\psi(x,0)=d(x,d(x,0))=d(x,x)=0$. We 
 then get  the following situation:
$$\xymatrix@=7pt{
	(x,y) \ar[dd]^{R[p_1]} \ar[rr]^{R[p_0]}  \ar@(l,l)@{.>}[dd]_{R[\psi]} && (x,0)\ar[dd]_{R[p_1]} \ar@(r,r)[dd]^{R[\psi]}\\
	&&&\\
	(0,y)  \ar[rr]_{R[p_0]} && (0,0)
}
$$
Whence  $d(x,d(x,y))=\psi(x,y)=\psi(0,y)=y$.

So, when $(X,d)$ is a subtraction, it is right cancellable as being the difference mapping of a group structure. Now, consider the following punctual span:
$$\xymatrix@=4pt{
	& X\times X  \ar@<-1ex>[ddd]_{p_1^X} \ar@<1ex>[rrrr]^{d}  &&&&  X  \ar@<-1ex>[ddd]_{\tau_X} \ar[llll]^{(1_X,0)} \\
	&&&&& && \\
	&&& &&&&&&\\
	& X  \ar[rrrr]^{\tau_X}\ar[uuu]_{s_0^X} &&&& 1 \ar@<1ex>[llll]^{0_X} \ar[uuu]_{0_X}
}
$$
It is a punctual relation since $d$ is right cancellable. Now we get the following situation, when $d'$ is another subtraction:
$$\xymatrix@=7pt{
	(x,y) \ar[dd]^{R[p_1]} \ar[rr]^{R[d]}  \ar@(u,u)@{.>}[rr]_{R[d']} && (d(x,y),0)\ar[dd]_{R[p_1]} \\
	&&&\\
	(y,y) \ar@(d,d)[rr]^{R[d']} \ar[rr]_{R[d]} && (0,0)
}
$$
whence the upper bended arrow which means $d'(x,y)=d(x,y)$.

Now, let $(Y,d)$ another subtraction and $f: X\to Y$ any map in $\EE$. Then consider the following diagram:
$$\xymatrix@=4pt{
	X\ar[ddd]_{f} \ar@{ >->}[rrrr]^{s_0^X}	&&&& X\times X  \ar[ddd]_{f\times f} \ar[rrrr]^{d}  &&&&  X   \ar@{.>}[ddd]^{\phi}   \\
	&&&&& && \\
	&&& &&&&&&\\
	Y\ar@{ >->}[rrrr]_{s_0^Y}	&&&& Y\times Y  \ar[rrrr]_{d} &&&& Y  
}
$$
Since we have horizontal cokernel diagrams, we get a unique  factorization $\phi$, satisfying $\phi(d(a,b))=d(f(a),f(b))$. Whence $\phi(a)=\phi(d(a,0))=d(f(a),0)=f(a)$. Finally we have: $x\circ x=0 \iff d(x,d(0,x))=0 \iff x=d(0,x)$.
\endproof

The uniqueness property asserted by this proposition makes intrinsic the notion of internal abelian group object in a congruence hyperextensive category $\EE$, whence the following:
\begin{defi}
	Let $\EE$ be a congruence hyperextensible category. An object  $X$ is to  be abelian when it is endowed with a (necessarily unique) substraction $d$. We shall denote $Ab\EE$ the  full subcategory of abelian  objects in $\EE$ which is additive and stable under finite limits in $\EE$. 
\end{defi}
When $\VV$ is a congruence hyperextensible variety, the subcategory $Ab\VV$ is actually abelian. Later on, we shall need the following observation:
\begin{prop}\label{birkh}
Let $\EE$ be a regular pointed category such that the category $Ab\EE$ of abelian groups in $\EE$ is a full subcategory of $\EE$. Suppose the inclusion $Ab\EE \into \EE$ has a left adjoint $A$. Then the following conditions are equivalent:\\
1) the universal comparison $\eta_X: X\to A(X)$ is a regular epimorphism;\\
2) $Ab\EE$ is stable under monomorphism;\\
3) $Ab\EE$ is Birkhoff subcategory of $\EE$. 	
\end{prop}
\proof
$\EE$ being finitely complete, $Ab\EE$ is stable under product. Suppose 1). Let $A$ be an abelian object and $u:U\into A$ a monomorphism. Then  the group homomorphism $A(u): A(U)\to A$ is such that $u=A(u).\eta_U$. The morphism $u$ being a monomorphism, so is $\eta_U$, which, being a regular epimorphism as well, is an isomorphism. Accordingly, $U$ is abelian. Suppose  2). Take the canonical decomposition of $\eta_X$, it produces a subobject $v: V\into A(X)$, which according to 2) is a subgroup of the abelian group $A(X)$. The  universal property of $A(X)$ makes $v$ an isomorphism and  $\eta_X$ a regular epimorphism.
 
We are now going to show that under the condition 2), any regular epimorphism $f: A\onto Y$ with $A$ abelian makes $Y$ abelian. Under the Condition 2), the object $R[f]$, as a subobject of $A\times A$, is abelian. Consider now the following diagram:
$$\xymatrix@=10pt{
	R[f]\times R[f] \ar@<-6pt>[rrr]_{d_0^f\times d_0^f} \ar@<6pt>[rrr]^{d_1^f\times d_1^f}
	\ar@<-8pt>[dd]_{p_0}  \ar[dd]|{\circ} \ar@<8pt>[dd]^{p_1}	&&& A\times A\ar@{ >->}[lll] \ar@<-8pt>[dd]_{p_0}   \ar@<8pt>[dd]^{p_1} \ar[dd]|{\circ} \ar@{->>}[rrrr]^{f\times f}  &&&& Y\times Y\ar@<6pt>[dd]^{p_1} \ar@<-6pt>[dd]_{p_0} \ar@{.>}[dd]|{\circ}\\
	&&&&&&&\\
	R[f] \ar@<-6pt>[rrr]_{d_0^f} \ar@<6pt>[rrr]^{d_1^f}
	   \ar@<6pt>[dd]^{\tau_{R[f]}}	&&& A \ar@{ >->}[lll]    \ar@<6pt>[dd]^{\tau_A} \ar@{->>}[rrrr]^{f}  &&&& Y \ar@<6pt>[dd]^{\tau_Y} \\
	&&&&&&&\\
	1 \ar@{=}[rrr]    \ar@{ >->}[uu]^{0_{R[f]}}	&&& 1  \ar@{=}[rrrr] \ar@{ >->}[uu]^{0_{A}} &&&& 1 \ar@{ >->}[uu]^{0_Y}
}
$$
In a regular category, the regular epimorphisms being stable under product, the morphism $f\times f$ is a regular epimorphism which, therefore, is the quotient of the upper horizontal left hand side equivalence relation. Accordingly the composition $\circ$ extends to the quotient and makes $Y$ an abelian group.
\endproof

\begin{prop}\label{comob}
	Let $\EE$ be a congruence hyperextensive category. Any binary operation $\circ$ with unit on an object $X$ is necessarily associative.
\end{prop}
\proof
Let $\circ$ be any binary operation with unit on  $X$. Consider the following punctual relation:
$$\xymatrix@=4pt{
	& X\times X \times X  \ar@<-1ex>[ddd]_{p_0^X} \ar@<1ex>[rrrrr]^>>>>>>>>{p_1^X.p_2^X}  &&&&&  X  \ar@<-1ex>[ddd]_{\tau_X} \ar[lllll]^<<<<<<<{0_{X\times X}\times X} \\
	&&&&& && \\
	&&& &&&&&&\\
	& X\times X  \ar[rrrrr]^>>>>>>>>>>{\tau_{X\times X}}\ar[uuu]_{X\times X \times 0_X} &&&&& 0 \ar[uuu]_{0_X} \ar@<1ex>[lllll]^<<<<<<<<<<{0_{X\times X}}
}
$$
Setting $\gamma(x,y,z)=x+(y+z)$, we get the following situation:
$$\xymatrix@=7pt{
	(x,y,z) \ar[dd]^{R[p_1.p_2]} \ar[rr]^{R[p_0]}  \ar@(l,l)@{.>}[dd]_{R[\gamma]} && (x,y,0)\ar[dd]_{R[p_1.p_2]} \ar@(r,r)[dd]^{R[\gamma]}\\
	&&&\\
	(0,x+y,z)  \ar[rr]_{R[p_0]} && (0,x+y,0)
}
$$
When $\EE$ is congruence  hyperextensive, we get: $\gamma(x,y,z)=\gamma(0,x+y,z)$, namely $x+(y+z)=(x+y)+z$.
\endproof

\subsection{Congruence hypoextensible categories}

\begin{defi}
	A pointed category $\EE$ is said to be congruence hypoextensible when given any punctual span in $\EE$ and any equivalence relation $T$ on $W$ such that $R[f]\cap R[g]\subset T\subset R[f]$, we get: $T \subset g^{-1}(t^{-1}(T))$.  
\end{defi}
In other words, when the following situation holds:

$$ \xymatrix@=15pt{
	x \ar[r]^{R[g]} \ar@(l,l)[d]_{T} \ar[d]_{R[f]} & tg(x) \ar[d]^{R[f]}  \ar@(r,r)@{.>}[d]^{T}\\
	x' \ar[r]_{R[g]} & tg(x') 
}
$$
Immediately we get the following remarkable categorical  property:
\begin{prop}\label{dimhypo}
Let $\EE$ be a congruence hypoextensible category. Given any punctual span in $\EE$, then any equivalence relation $T$ on $W$ such that $R[f]\cap R[g]\subset T\subset R[f]$ has a direct image along $g$ which is nothing but $t^{-1}(T)$.
\end{prop}
\proof
The previous diagram provides us with the dotted morphism $\check g$ of equivalence relations:
$$\xymatrix@=2pt{
	T \ar@<-2ex>[dddd]_{d_0^T}     \ar@<2ex>[dddd]^{d_1^T} \ar@<1ex>@{.>}[rrrrr]^{\check g}  &&&&&  t^{-1}(T)  \ar@<-2ex>[dddd]_{d_0} \ar[lllll]^{\check t}\ar@<2ex>[dddd]^{d_1} \\
	&&&&& &&&&\\
	&&& &&&&&&\\
	&&& &&&&&&\\
	W \ar[uuuu]|{s_0^T}  \ar@<1ex>[rrrrr]^{g}  &&&&& Z \ar[lllll]^{t} \ar[uuuu]|{s_0} 
}
$$
defined by $\check g(uTv)=(g(u),g(v))$, which is regular epimorphism in $\EE$ since it is split by $\check t$.
\endproof

Any pointed subtractive variety $\VV$ \cite{Ur} or regular subtractive category $\EE$ is congruence hypoextensible \cite{BGum}, where the concept of subtractive category, introduced in \cite{ZJ}, is a way to categorically characterize the pointed subtractive varieties:
\begin{defi}
A pointed category $\EE$ is subtractive when any split left punctual relation is a punctual one.
\end{defi}
Here is an example strictly weaker than the pointed substractive varietal context which is  obtained when we have $\epsilon=Id$:
\begin{prop}\label{exhypo}
Let $\VV$ be a pointed variety with a binary term $d$ and a unary term $\epsilon$ such that:
\begin{align*}
d(a,a)=0\;\;\;\;\; ;\;\;\;\;  \;\;\; \epsilon(d(a,0)=a 
\end{align*}
Then the variety $\VV$ is congruence hypoextensible. 
\end{prop}
\proof
Consider the following punctual relation in $\VV$:
$$\xymatrix@=4pt{
& W  \ar@<-1ex>[ddd]_{\pi_X} \ar@<1ex>[rrrr]^{\pi_Z}  &&&&  Z  \ar@<-1ex>[ddd]_{\tau_Z} \ar[llll]^{0_Z^W} \\
&&&&& && \\
&&& &&&&&&\\
& X  \ar[rrrr]^{\tau_X}\ar[uuu]_{0_X^W} &&&& 1 \ar@<1ex>[llll]^{0_X} \ar[uuu]_{0_Z}
}
$$
and an equivalence relation $T$ on $W$ such that $(aWb)T(aWc)$. Then we get: 
\begin{align*}
d(aWb,aW0)&=d(a,a)Wd(b,0)=d(0,0)Wd(b,0)&=d(0Wb,0W0) \\
\; \;{\rm whence:} \;\; & \;\;\;\;\;\;\;\;\; d(0Wb,0W0)Td(0Wc,0W0) &
\end{align*}
Then, applying $\epsilon$ to each side, we get: $(0Wb)T(0Wc)$.
\endproof

\begin{prop}\label{2}
	When $\EE$ is congruence hypoextensible, then any internal binary operation $\circ$ on an object $X$ is left cancellable as soon as it has a left unit. In particular any opsubtraction is left cancellable, and, by duality, any subtraction is right cancellable.
\end{prop}
\proof
We get the following situation:
$$\xymatrix@=7pt{
	(a,c) \ar[dd]^{R[p_0]} \ar[rr]^{R[p_1]}  \ar@(l,l)[dd]_{R[\circ]} && (0,c)\ar[dd]_{R[p_0]} \ar@(r,r)@{.>}[dd]^{R[\circ]}\\
	&&&\\
	(a,c')  \ar[rr]_{R[p_1]} && (0,c')
}
$$
which means the left cancellability: $a\circ c=a\circ c' \Rightarrow c=c'$. 
\endproof

\subsection{Punctually congruence modular categories}

\begin{defi}
	A punctually congruence modular category $\EE$ is a pointed category which is both congruence hyperextensible and hypoextensible; i.e. such that, given any punctual span and any equivalence relation $T$ on $W$ such that $R[f]\cap R[g]\subset T\subset R[f]$, we get: $T=R[f]\cap g^{-1}(t^{-1}(T))$.
\end{defi}
The categorical translation is the following:
\begin{prop}
	Let $\EE$ be a punctually congruence modular category. Given a punctual span in $\EE$ and an equivalence relation $T$ on $W$  such that $R[f]\cap R[g]\subset T\subset R[f]$, the following diagram where $\check g$ is defined by Proposition \ref{dimhypo}:
	$$\xymatrix@=1pt{
		T \ar@{ >->}[dddd]   \ar@<1ex>[rrrrr]^{\check g}  &&&&&  t^{-1}(T)  \ar@{ >->}[dddd]_{} \ar[lllll]^{\check t} \\
		&&&&& &&&&\\
		&&& &&&&&&\\
		&&& &&&&&&\\
		R[f]   \ar@<1ex>[rrrrr]^{R(g)}  &&&&& \nabla_Z \ar[lllll]^{R(t)}  
	}
	$$
   a is pullback in the category $Equ\EE$ of equivalence relations in $\EE$.
\end{prop}

Any strongly unital variety (= subtractive J\'onsson-Tarsky variety) $\VV$ or any regular strongly unital (=unital+ subtractive, see \cite{B0} and \cite{ZJ})) category $\EE$ is punctually congruence modular \cite{BGum}. The following example is stricly weaker than the strongly unital context which is obtained when we have $\upsilon=d(-,0)$ (and consequently $\epsilon=-\circ 0$):
\begin{prop}
	Let $\VV$ be a pointed variety with two binary terms $(\circ,d)$ and two unary terms $\upsilon$ and $\epsilon$ such that:
	\begin{align*}
	a\circ0=0\circ a; \;\; \;\;\;\;\; \;\; d(a,a)=0; &\;\;\;\;\;\;\;\;  \upsilon(a\circ 0)=a; \;\;\;\;\;  \;\;\;\;\; \epsilon(d(a,0))=a
	\end{align*}
	Then the variety $\VV$ is punctually congruence modular.
\end{prop}
\proof
Straightforward from Propositions \ref{exhyper} and \ref{exhypo}.
\endproof

\begin{prop}\label{cancel}
Let $\EE$ be a punctually congruence category. 
A binary operation $(X,\circ)$ with a right hand side unit has a unique unary right hand side canceller $!: X\to X$, namely a unique unary  operation such that $x\circ !x=0$.
	
An object $X$ is abelian in $\EE$ if and only if it is endowed with a binary operation $\circ$ which has a right hand side unit and a unary right canceller.
\end{prop}
\proof
Any binary operation $(X,\circ)$ with a right hand side unit and a right hand side canceller $!: X\to X$ produces a subtraction $d(x,y)=x\circ !y$ on $X$. Since $d$ is unique, then $!(x)=d(0,x)$ is unique. This subtraction $d$ makes $X$ abelian.
\endproof

\section{The diagonal punctuation}

It is worth, now, dwelling on the specific example of a congruence hypoextensive categories given by the regular subtractive categories and, accordingly, by the pointed subtractive varieties, since, in this context, we have two meaningful major properties which do not seems to hold in the mere case of regular congruence hypoextensive categories, but which will appear to be valid in any pointed fiber $Pt_Y\VV$ of a congruence modular  variety.
\begin{prop}[\cite{BZZ}]
	In a subtractive category $\EE$, on any object $X$, there is at most one abelian group structure, and the subcategory $Ab\EE$ of abelian objects in $\EE$ is full and stable under monomorphisms.
\end{prop} 
\begin{prop}[\cite{BZ}]
	In a regular subtractive category $\EE$, the abelianization of any object $X$, when it exists, is given by the cokernel of the diagonal $s_0^X:X\to X\times X$.
\end{prop}
So, let us now introduce the following guide-line tool:
\begin{defi}
	Let $X$ be an object in any category $\EE$. Its diagonal punctuation $DpX$ is defined, when it exists, by the following left hand side upward pushout:
	$$\xymatrix@=7pt{
		X\times X   \ar[rr]^{\omega_X}  && DpX &&& X\times X \ar@<-6pt>[dd]_{p_0^X}  \ar[rr]^{\omega_X}  && DpX \ar@<6pt>[dd]_{}\\
		&&&&&&\\
		X \ar@{ >->}[uu]^{s_0^X} \ar[rr]_{\tau_X} && 1 \ar@{ >->}[uu]_{0} &&& 	X \ar@{ >->}[uu]_{s_0^X}  \ar[rr]_{\tau_X} && 1 \ar@{ >->}[uu]^{0}
	}
	$$
\end{defi}
Then the retraction $p_0^X$ produces the downward pushout on the right hand  side.  Of course, we have a similar downward pushout produced by the retraction $p_1^X$. The twisting isomorphism $tw_X(u,v)=(v,u)$ induces an  isomorphism $\gamma_X$ on $DpX$ such that $\gamma_X.\omega_X=\omega_X.tw_X$ and $\gamma_X(0)=0$. This morphism is involutive since so is $tw_X$. So, we get immediately:
\begin{lemma}\label{diag2}
	For any couple $(f,g),(f',g')$ of parallel pairs of morphisms with codomain $X$  we have: $\omega_X(f,g)=\omega_X(f',g')\iff \omega_X(g,f)=\omega_X(g',f')$.
\end{lemma}
\proof
We get: $\omega_X(g,f)=\gamma_X.\omega_X(f,g)=\gamma_X.\omega_X(f',g')=\omega_X(g',f')$.
\endproof

Let us now consider the previous definition diagram completed by the horizontal kernel equivalence relations:
$$\xymatrix@=10pt{
	R[\omega_X] \ar@<-6pt>[rrr]_{d_0^{\omega}} \ar@<6pt>[rrr]^{d_1^{\omega}}
	\ar@<-6pt>[dd]_{R(p_0)}   \ar@<6pt>[dd]^{R(p_1)}	&&& X\times X \ar@{ >->}[lll] \ar@<-6pt>[dd]_{p_0^X}   \ar@<6pt>[dd]^{p_1^X} \ar[rr]^{\omega_X}  && DpX \ar@<6pt>[dd]^{} \\
	&&&&&&\\
	X\times X \ar@<-6pt>[rrr]_{p_0^X}   \ar@<6pt>[rrr]^{p_1^X} \ar@{ >->}[uu]	&&& X \ar@{ >->}[uu]|{s_0^X} \ar[rr]_{\tau_X} \ar@{ >->}[lll] && 1 \ar@{ >->}[uu]|{0}
}
$$
Then this diagram produces the vertical left hand side reflexive relation $D_{\omega}$ on $X\times X$. In set-theoretical terms, we get $(x,y)D_{\omega}(x',y') \iff \omega_X(x,x')=\omega_X(y,y')$.
\begin{lemma}\label{diag3}
	The relation $D_{\omega}$ is symmetric, and we have: $(s_0^X)^{-1}(D_{\omega})=\nabla_X$.
\end{lemma}
\proof
The first point is a straighforward consequence of Lemma \ref{diag2}; and the second one of:  $D_{\omega}(x,x)=D_{\omega}(y,y) \iff \omega_X(x,y)=\omega_X(x,y)$
\endproof

\begin{theo}\label{DD}
	Let $\EE$ be a Gumm category with diagonal punctuations. Then:\\
	1) we get: $\omega_{X}(t,t)=\omega_{X}(t,t')\iff  \omega_{X}(x,t)=\omega_{X}(x,t'), \forall x\in X$;\\
	2) the reflexive and symmetric relation $D_{\omega}$ is actually an  equivalence relation. This result holds, a fortiori, in any congruence modular variety.
\end{theo}
\proof
1) Consider the following diagram  in $X\times X$:
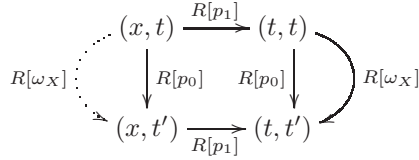
$$\xymatrix@=8pt{
	(x,t)  \ar[dd]^{R[p_0]}   \ar@{.>}@(l,l)[dd]_{R[\omega_X]}\ar[rr]^{R[p_1]} && (t,t) \ar@(r,r)[dd]^{R[\omega_X]}  \ar@<1ex>[dd]_{R[p_0]}\\
	&&&&\\
	(x,t')    \ar[rr]_{R[p_1]} && (t,t')    
}
$$

2) It remains to check that the relation $D_{\omega}$ is transitive. So suppose that $(x,y)D_{\omega}(x',y')D_{\omega}(x'',y'')$, namely that  $\omega(x,x')=\omega(y,y')$ and $\omega(x',x'')=\omega(y',y'')$.
We have to conclude that $\omega(x,x'')=\omega(y,y'')$. Using the symmetry of $D_{\omega}$, it is a consequence of the \emph{Cube  Lemma} \ref{cubl} applied to the following situation:
$$
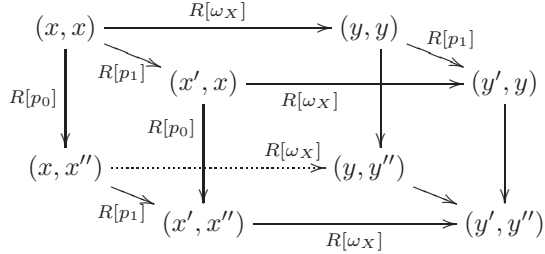
\xymatrix@=5pt{
	(x,x) \ar[rrd]_{R[p_1]}  \ar[ddd]_{R[p_0]}   \ar[rrrrr]^{R[\omega_X]} &&&&& (y,y) \ar[rrd]^{R[p_1]}  \ar@<1ex>[ddd]^{}\\
	&& (x',x) \ar[ddd]_<<<<{R[p_0]}  \ar[rrrrr]_<<<<<<<<<{R[\omega_X]}  &&&&& (y',y) \ar[ddd]^{}  && \\
	&&& &&&&&&\\ 
	(x,x'') \ar[rrd]_{R[p_1]}  \ar@{.>}[rrrrr]^>>>>>{R[\omega_X]}   &&&&& (y,y'')  \ar[rrd]  \\
	&&  (x',x'')  \ar[rrrrr]_{R[\omega_X]}    &&&&& (y',y'') 
}
$$
since it produces the desired dotted arrow. 
\endproof 

\begin{lemma}\label{diag4}
	Let $\EE$ be a Gumm category with diagonal punctuations. The diagonal $s_0^X$ lies inside an equivalence class of $D_{\omega}$.
\end{lemma}
\proof
It is a straightforward consequence of the previous lemma.
\endproof
Let us now come back to the varietal context. Any congruence modular variety $\VV$ is factor permutable \cite{Gu}. From that we shall deduce two important consequences:
\begin{prop}\label{imp1}
	Let $\VV$ be a congruence modular variety. Then for any non-empty algebra $X$, the diagonal punctuation diagram produces a downward \emph{regular} pushout:
	$$\xymatrix@=7pt{
		X\times X  \ar@<-6pt>[dd]_{p_0^X}    \ar@{->>}[rr]^{\omega_X}  && DpX \ar@<6pt>[dd]^{} \\
		&&&&\\
		X \ar@{ >->}[uu]|{s_0^X} \ar@{->>}[rr]_{\tau_X}  && 1 \ar@{ >->}[uu]|{0}
	}
	$$	
\end{prop}
\proof
We have to show that the factorization $(p_0^X,\omega_X): X\times X\to X\times DpX$ is a surjective homomorphism. So, from a triple $(x,a,b)$ of elements of $X$, we have to find an element $y$ such that $\omega_X(x,y)=\omega_X(a,b)$. For that consider the following situation depicted by the plain arrows:
$$\xymatrix@=5pt{
	(x,x) \ar[dd]_{R[\omega_X]} \ar@{.>}[rr]^{R[p_0^X]}   &&  (x,y) \ar@{.>}[dd]^{R[\omega_X]}\\
	&&&\\
	(a,a) \ar[rr]_{R[p_0^X]} && (a,b)
}
$$
The permutation of the pair $(R[p_0^X],R[\omega_X])$ produces the dotted part of the diagram which gives us the $y$ in question.
\endproof

\begin{prop}\label{imp2}
	Let $\VV$ be a congruence modular variety. Then for any non-empty algebra $X$, the pointed algebra $(DpX,0)$ is endowed with an internal group structure in $\VV$ which is necessarily abelian and the unique one in $\VV$ with unit $0$. 
\end{prop}
\proof
Consider now the previous diagram completed, at the upper level, by the second level of the two  vertical left hand side equivalence relations. This produces the upper horizontal left hand side equivalence relation. 
$$\xymatrix@=8pt{
	R[R[p_0]] \ar@<-6pt>[rrr]_{R(d_0^{\omega})} \ar@<6pt>[rrr]^{R(d_1^{\omega})}
	\ar@{=}[dd]   	&&& R[p_0^X] \ar@{ >->}[lll] \ar@{=}[dd] \ar@{-->}[rrr]^{(\omega_X.p_0^X,\omega_X.p_1^X)}  &&& DpX \times DpX\ar@{=}[dd]\\
	&&&&&&\\
	(D_{\omega})_2 
	\ar@<-8pt>[dd]_{d_0}  \ar[dd]|{d_1} \ar@<8pt>[dd]^{d_2}	&&& X\times X\times X \ar@<-8pt>[dd]_{p_0^X}   \ar@<8pt>[dd]^{p_2^X} \ar[dd]|{p_1^X}   &&& DpX \times DpX \ar@<6pt>@{.>}@(r,r)[dd]^{\partial_2} \ar@<-6pt>[dd]_{p_0} \ar[dd]^{p_1}\\
	&&&&&&\\
	R[\omega_X] \ar@<-6pt>[rrr]_{d_0^{\omega}} \ar@<6pt>[rrr]^{d_1^{\omega}}
	\ar@<-6pt>[dd]_{R(p_0)}   \ar@<6pt>[dd]^{R(p_1)}	&&& X\times X \ar@{ >->}[lll] \ar@<-6pt>[dd]_{p_0^X}   \ar@<6pt>[dd]^{p_1^X} \ar@{->>}[rrr]^{\omega_X}  &&& DpX \ar@<6pt>[dd]^{} \\
	&&&&&&\\
	X\times X \ar@<-6pt>[rrr]_{p_0^X}   \ar@<6pt>[rrr]^{p_1^X} \ar@{ >->}[uu]	&&& X \ar@{ >->}[uu]|{s_0^X} \ar@{->>}[rrr]_{\tau_X} \ar@{ >->}[lll] &&& 1 \ar@{ >->}[uu]|{0}
}
$$
The upper horizontal right hand side map is the one which makes commute the upper right hand side squares indexed by $0$ and $1$. So, by commutation of limits, the upper horizontal diagram produces an effective equivalence relation, namely a map completed with its kernel equivalence relation. Now, according to the previous proposition this same map $R[p_0^X]\to DpX\times DpX$, induced by the lower right hand side regular pushout, is a regular epimorphism; namely it is the quotient of its kernel equivalence relation. So the left hand side squares indexed by $2$ induces the vertical right hand side $\partial_2$ which gives a  group structure with  unit $0$ to the object $DpX$.

This group structure is necessarily abelian and the unique one with unit $0$, according to Proposition \ref{abob} applied to the congruence hyperextensive fiber $Pt_1\EE$. So, we get:
$$ (*) \; -\omega_X(a,b)+\omega_X(a,c)=\omega_X(b,c);\; {\rm whence:\;} \omega_X(a,c)=\omega_X(a,b)+\omega_X(b,c)$$

The diagonal punctuation determines the following property:
if $(A,\circ,0)$ is another internal abelian group in $\VV$, given any $\VV$-homomorphism  $f:X\to A$, we get the following large commutative diagram in $\VV$ where $\partial_A$ is the  \emph{opdifference mapping} of the group $A$:
$$\xymatrix@=5pt{
	& A\times A \ar[rrr]^{\partial_A} &&& A\\
	&&&&\\
	X\times X   \ar@{->>}[rr]^{\omega_X} \ar[ruu]^{f\times f}  && DpX \ar@{.>}[rruu]^{\bar f} \\
	&&&&\\
	X \ar@{ >->}[uu]^{s_0^X} \ar@{->>}[rr]_{\tau_X} && 1 \ar@{ >->}[uu]^{0}  \ar@<-1ex>@{ >->}[uuuurr]_{0}
}
$$
Since $\partial_A.s_0^A=0$, it produces a unique factorization $\bar f$ making the adjacent diagrams commute. Namely: $(1):\;\bar f\omega_X(a,b)=f(b)-f(a)$ and $(2):\;\bar f(0)=0$. Then, by Proposition \ref{abob}, the equality $(2)$ implies that $\bar f$ is certainly a group  homomorphism since, the map $\bar f$ is lying in the congruence hyperextensive  fiber $Pt_1\EE$ between two abelian objects, it necessarily preserves the abelian group structures.
\endproof
\begin{coro}
	Let $\VV$ be a congruence modular variety. 
	Given any non-empty algebra $X$ and its diagonal punctuation, the following conditions are equivalent:\\
	1) $s_0^X$ is the "kernel" of $\omega_X$, namely we have: $\omega_X(u,v)=0 \iff u=v$;\\
	2) the diagonal $s_0^X$ is a congruence class of $R[\omega_X]$;\\
	3) the following downward square is a pullback:
	$$\xymatrix@=8pt{
		X\times X  \ar@<-6pt>[dd]_{p_0^X}    \ar[rr]^{\omega_X}  && DpX \ar@<6pt>[dd]^{} \\
		&&&&\\
		X \ar@{ >.>}[uu]|{s_0^X} \ar[rr]_{\tau_X}  && 1 \ar@{ >.>}[uu]|{0}
	}
	$$
	4) $X$ is affine in $\VV$.  
\end{coro}
\proof
From Lemma \ref{diag4}, the equivalence $(1)\iff (2)$ is straightforward.
Since the above downward square is a regular pushout, it is a pullback if and only if the morphism $(p_0^X,\omega_X)$ is a monomorphism. By Theorem \ref{DD} we then get:
$\omega_{X}(t,t)=\omega_{X}(t,t')\iff  \omega_{X}(x,t)=\omega_{X}(x,t'), \forall x\in X$;
whence: $0=\omega_{X}(t,t')\iff  \omega_{X}(x,t)=\omega_{X}(x,t')$.
So, we get $(1)\iff (3)$.
Now, when we have (3), the comparison  $(p_0^X,\omega_X):X\times X\to X\times DpX$ is an isomorphism; let us denote $(p_0^X,\chi)$ its inverse. Then the ternary operation $p(u,v,w)=\chi(u,\omega_X(v,w))$ gives $X$ a Mal'tsev operation which insures that $X$ is affine.

Conversely, suppose that $X$ is affine. Then consider the following diagram, where $\vec A_X$ is the \emph{direction} of the affine object  $X$ (see \cite{B3} and Proposition \ref{abab} below), and consequently the downward quadrangle with $\vec A_X$ is a pullback:
$$\xymatrix@=8pt{
	&& && \vec A_X \ar@<12pt>[dddll]^{}\\
	X\times X  \ar@<-6pt>[dd]_{p_0^X}    \ar[rr]_{\omega_X}  \ar[rrrru]^{q_X} && DpX \ar@<6pt>[dd]^{} \ar@{.>}[rru]_{\epsilon} \\
	&&&&\\
	X \ar@{ >.>}[uu]|{s_0^X} \ar[rr]_{\tau_X}  && 1 \ar@{ >.>}[uu]|{0} \ar@{ >.>}@<-6pt>[uuurr]|{0}
}
$$
We get a factorization $\epsilon$ making the two triangles commute: $\epsilon.\omega_X=q_X$ and $\epsilon.0=0$. The first equality implies that $\epsilon$ is a regular epimorphism, since so is $q_X$. This $\epsilon$ determines a regular epimorphic factorization $\bar{\epsilon}:X\times DpX\to X\times X$ such that $q_X.\bar{\epsilon}=\epsilon.p_1$ and $p_0^X.\bar{\epsilon}=p_0$. So, $(p_0^X,\omega_X).\bar{\epsilon}=1_{X\times X}$ (by composition with the pair $(p_0^X,q_X)$), and $\bar{\epsilon}$ is a monomorphism as well. Consequently it is an isomorphism, and so is
$(p_0^X,\omega_X)$. Whence $(3)$. 
\endproof
So, in the pointed case, we recover a result of \cite{Gu}, but with a categorical proof:
\begin{coro}
	Let $\VV$ be a pointed congruence modular variety. A non-empty algebra $X$ is affine in $V$ if and only if its diagonal is a congruence class.
\end{coro}
\proof
Let $\Sigma$ be a congruence on $X\times X$ such that $s_0^X$ is the class of $(0,0)$ and $q:X\onto Q$ its quotient. Then, by the Barr-Kock Theorem, we get a pullback on the right hand side since we have pullbacks on the left hand side:
$$\xymatrix@=10pt{
	\Sigma \ar@<-6pt>[rrr]_{d_0^{\Sigma}} \ar@<6pt>[rrr]^{d_1^{\Sigma}}
	&&& X\times X \ar@{ >->}[lll]    \ar@{->>}[rrr]^{q}  &&& Q  \\
	&&&&&&\\
	X\times X  \ar@<-6pt>[rrr]_{p_0^X} \ar@<6pt>[rrr]^{p_1^X} \ar@{ >->}[uu]	&&& X \ar@{ >->}[uu]^{s_0^X} \ar@{->>}[rrr]_{\tau_X} \ar@{ >->}[lll] &&& 1 \ar@{ >->}[uu]_{0}
}
$$
Accordingly there is a unique factorization $\gamma: DpX\to Q$ such that $q=\gamma.\omega_X$, so that we get $R[\omega_X]\subset R[q]=\Sigma$. And,  since $(s_0^X)^{-1}(R[\omega_X])=\nabla_X$, the diagonal $s_0^X$ is a congruence class of $R[\omega_X]$ as well, which means that $\omega_X(x,y)=0$ if and only if $x=y$ or, in other words, means that $s_0^X$ is the "kernel" of $\omega_X$.
\endproof
In the pointed case, we get the following precision:
\begin{coro}
	When $\VV$ is a pointed congruence modular variety, $DpX$ is the universal abelian group associated with $X$ in $\VV$, and the universal comparison $\eta_X: X\onto DpX$ is a surjective homomorphim. Accordingly the subvariety $Ab\VV \into \VV$ is a Birkhoff one.
\end{coro}
\proof 
Let us show that $\omega_X.(0,1_X):X\to DpX$ is the universal comparison with the abelian groups in $\VV$. So, let $f: X\to A$ be any $\VV$-homomorphism, and consider the map $\bar f :X\times X \to A$ described in Proposition \ref{imp2}. From $(1)$, we get $\bar f\omega_X(0,a)=f(a)$ which means: $(**): f=\bar f.\omega_X(0,1_X)$. It remains  to show the uniqueness of such a factorization satisfying $(**)$. So, let $\phi: DpX\to A$ be another group homomorphism such that $f=\phi.\omega_X(0,1_X)$. Then:
\begin{align*}
\phi(\omega_X(a,b))=phi(\omega_X(0,b)-\omega_X(0,a))=&\;\phi(\omega_X(0,b))-\phi(\omega_X(0,a))\\
=f(b)-f(a)=\bar f(\omega_X(a,b))\; & 
\end{align*}
Now, since $\omega_X$ is a surjective homomorphism, we finally get the desired uniqueness $\phi=\bar f$. Finally to show that the universal  comparison $\omega_X.(0,1_X): X\to DpX$ is sujective, apply Proposition \ref{imp2}  to the triple $(0,a,b)$. The  last sentence is then a consequence of Proposition \ref{birkh}.
\endproof
Any pointed $3$-permutable variety is an example of pointed congruence modular variety; it is the case in particular of the variety $IA$ of pointed idempotent implication algebras. The previous precision can be now transfered to the non-pointed case. Here is now the result we were aiming to:

\begin{theo}\label{main}
	Let $\VV$ be a congruence modular variety. Given any algebra $Y$ and any split epimorphism $(f,s):X\splito Y$, consider the diagonal punctuation construction in the fiber $Pt_Y\VV$, namely the pushout of $s_0^f$ along $f$:
	$$\xymatrix@=8pt{
		R[f]  \ar@<-6pt>@{.>}[dd]_{d_0^f}  \ar@{->>}[rrr]^{\omega_f}  &&& Dp[f] \ar@<6pt>@{.>}[dd]^{\psi_f} \\
		&&&&&&\\
		X \ar@{ >->}[uu]_{s_0^f} \ar@{->>}[rrr]_{f}  &&& Y \ar@{ >->}[uu]^{\theta_f}
	}
	$$
	and let $\psi_f$ denote the retraction induced by $d_0^f$. Then:\\
	1) the downward square is a \emph{regular} pushout;\\
	2) the split epimorphism $(\psi_f,\theta_f)$ is the universal abelian split epimorphism associated with $(f,s)$.\\
	The subcategory $APt_Y\VV$ of  abelian split epimorphisms above $Y$ is  stable under monomorphism in $Pt_Y\EE$.
\end{theo}
\proof
Since $(f,s)$ is a split epimorphism in $\VV$, the homomorphism $f$ is a regular epimorphism, and the pushout construction makes $\omega_f$ a regular epimorphism. We have to show that the map $(d_0^f,\omega_f): R[f]\to X\times_YDp[f]$ is a regular epimorphism, namely a surjective one.  For that, consider the following pullback of split epimorphisms in $\VV$:
$$\xymatrix@=8pt{
	R[f]    \ar[dd]_{d_0^f} \ar@{->>}[rrrr]^{d_1^f}  &&&& X \ar[dd]_{f} \\
	&&&&&&&\\
	X \ar@{ >->}@<-6pt>[uu]_{(1_X,sf)} \ar@{->>}[rrrr]_{f}  &&&& Y \ar@<-6pt>@{ >->}[uu]_{s}
}
$$
Then, since it is a pullback,  we get the permutation of the pair $(R[d_0^f],R[d_1^f])$. Since it is a commutative diagram of split epimorphims, we have in addition: $R[d_0^f]\vee R[d_1^f]=R[f.d_0^f]$. Now the identity  $\psi_f.\omega_{f}=f.d_0^f=f.d_1^f$ shows that we get $R[\omega_{f}]\subset R[f.d_0^f]=R[d_0^f]\vee R[d_1^f]$. So according to Corollary 4.4 in \cite{Gu}, the pair $(R[\omega_{f}],R[d_0^f])$ does permute.

Now, starting with a triple $(x,a,b)$ of elements in $X$ such that $f(x)=f(a)=f(b)$, we have to find an element $y$ in $X$ such that $f(y)=f(x)$ and $\omega_{f}(x,y)=\omega_{f}(a,b)$. By the above permutation we get the following situation:
$$\xymatrix@=5pt{
	(x,x) \ar[dd]_{R[\omega_f]} \ar@{.>}[rr]^{R[d_0^f]}   &&  (x,y) \ar@{.>}[dd]^{R[\omega_f]}\\
	&&&\\
	(a,a) \ar[rr]_{R[d_0^f]} && (a,b)
}
$$
which produces the desired $y$ in question. 

The end of the proof is now similar to the end of the proof of Proposition \ref{imp2}, but, now, inside the slice category $\VV/Y$. First complete the diagram of the diagonal punctuation by  the horizontal kernel equivalence relations:
$$\xymatrix@=10pt{
	(D_{\omega_f})_2 \ar@<-6pt>[rrr]_{R(d_0^{\omega})} \ar@<6pt>[rrr]^{R(d_1^{\omega})}
	\ar@<-8pt>[dd]_{d_0}  \ar[dd]|{d_1} \ar@<8pt>[dd]^{d_2}	&&& R^2[f]\ar@{ >->}[lll] \ar@<-8pt>[dd]_{d_0^f}   \ar@<8pt>[dd]^{d_2^f} \ar[dd]|{d_1^f} \ar@{-->}[rrrr]^{(\omega_{(f,s)}.d_0^f,\omega_{(f,s)}.d_1^f)}  &&&& R[\psi_f] \ar@<6pt>@{.>}@(r,r)[dd]^{\partial_2} \ar@<-6pt>[dd]_{d_0} \ar[dd]^{d_1}\\
	&&&&&&&\\
	R[\omega_{f}] \ar@<-6pt>[rrr]_{d_0^{\omega}} \ar@<6pt>[rrr]^{d_1^{\omega}}
	\ar@<-6pt>[dd]_{R(d_0^f)}   \ar@<6pt>[dd]^{R(d_1^f)}	&&& R[f] \ar@{ >->}[lll] \ar@<-6pt>[dd]_{d_0^f}   \ar@<6pt>[dd]^{d_1^f} \ar@{->>}[rrrr]^{\omega_f}  &&&& Dp[f] \ar@<6pt>@{.>}[dd]^{\psi_f} \\
	&&&&&&&\\
	R[f] \ar@<-6pt>[rrr]_{d_0^f}   \ar@<6pt>[rrr]^{d_1^f} \ar@{ >->}[uu]	&&& X \ar@{ >->}[uu]|{s_0^f} \ar@{->>}[rrrr]_{f} \ar@{ >->}[lll] &&&& Y \ar@{ >->}[uu]^{\theta_f}
}
$$
The variety $\VV$ being congruence modular, the slice category $\VV/Y$ is a Gumm category, which is finitely cocomplete. So, by Proposition \ref{DD}, the vertical left hand side reflexive relation is actually an equivalence relation we shall denote by $D_{\omega_f}$. 

Then complete the diagram at the upper level by the second step of the two vertical equivalence relations on the left hand side. Now the map $\psi=(\omega_{f}.d_0^f,\omega_{f}.d_1^f):R^2[f]\to  R[\psi_{f}]$ is defined by: $$\psi(a,b,c)=(\omega_{f}(a,b),\omega_{f}(a,c))$$ It is the one which makes commute the pairs $(d_0^f,d_1^f)$ with the pairs  indexed by $0$ and $1$ on the vertical right hand side. This  map $\psi$ is a surjective homomorphism since the square below is a regular pushout. So, it determines a factorization  $\partial_2$ on the  right hand side, giving a group structure (with unit $\theta_f$) to the object  $\psi_{f}$ in the slice category $\VV/Y$.  This group structure is necessarily abelian and the unique one with unit $\theta_{f}$ in the Gumm category $\VV/Y$, acccording to Proposition. So, we get:
$ (*) \; -\omega_{f}(a,b)+\omega_{f}(a,c)=\omega_{f}(b,c)$,
and the following identities:
\begin{align*}
(\bar 1)\;\;\;\; \omega_X(a,b)+\omega_X(b,sf(b))=\omega_X(a,sf(b)) &=\omega_X(a,sf(a)) \\ 
(\bar 2)\;\; \omega_X(a,b)=\omega_X(a,sf(a))-\omega_X(b,sf(b));& \; (\bar 3)\; \omega_X(sf(b),b)=-\omega_X(b,sf(b)) 
\end{align*}  

Let us show now that the homomorphism $\omega_{f}.(sf,1_X): X\onto Dp(f,s)$ produces the universal comparison with the abelian split epimorphism $(\psi_f,\theta_f)$. First let us show it surjective. For that apply the previous reasonning to:
$$\xymatrix@=5pt{
	(x,y) \ar@{.>}[rr]^<<<<{R[\omega_{f}]} \ar[dd]_{R[d_0^f]}   &&  (sf(x),t) \ar@{.>}[dd]^{R[d_0^f]} \ar@{=}[r] & ((sf(t),t))\\
	&&&\\
	(x,x) \ar[rr]_<<<<{R[\omega_{f}]} && (sf(x),sf(x))
}
$$
Given any abelian split epimorphism $(g,t): A\splito Y$ and any $\VV$ homomorphism $h: X\to A$ between $(f,s)$ and $(g,t)$, let us consider the following diagram where $\partial:R[g]\to A$ is the opsubtraction which gives the abelian structure to the split epimorphism $(g,t)$:
$$\xymatrix@=5pt{
	& R[g] \ar[rrr]^{\partial} &&& A\\
	&&&&\\
	R[f]   \ar@{->>}[rr]^{\omega_{f}} \ar[ruu]^{R(h)}  && Dp[f] \ar@{.>}[rruu]^{\bar h} \\
	&&&&\\
	X \ar@{ >->}[uu]^{s_0^X} \ar@{->>}[rr]_{f} && Y \ar@{ >->}[uu]^{\theta_{f}}  \ar@<-1ex>@{ >->}[uuuurr]_{t}
}
$$
The pushout definition of $Dp[f]$ produces a  unique factorization $\bar h$, making commute the two adjacent diagrams. The commutation of the right hand side triangle makes $\bar h$ a group homomorphism in the congruence hyperextensible fiber $Pt_Y\VV$, while the other quadrangle means: $\bar h\omega_{(f,s)}(x,y)=\partial(h(x),h(y))$. Whence:
$\bar h\omega_{f}(sf(x),x)=\partial(hsf(x),h(x))=\partial(tgh(x),h(x))=h(x)$.
From that we get $g.\bar h=\psi_{f}$ and $\bar h.\theta_{f}=t$ by composition with the surjective homomorphism $\omega_{f}$. The uniqueness of this factorization comes from the fact that $\omega_{f}.(sf,1_X)$ is a surjective homomorphism as well. This same fact, by Proposition \ref{birkh},  insures that $APt_Y\VV$ is stable under monomorphisms in $Pt_Y\VV$.
\endproof

\section{Gumm categories}\label{kGumm}

\subsection{Direct consequence of the characterization theorem}

A first consequence of this theorem is to produce, thanks to Theorem \ref{teo}, a relatively simple varietal congruence modular setting from Proposition \ref{exhyper}:

\begin{prop}
	Let $\VV$ be a variety with a ternary term $p$ and a binary term $\alpha$ such that:
	\begin{align*}
	p(x,y,y)=p(y,y,x) &\;\;\;{\rm and}\;\;\;\;\;\; p(x,x,x)=x \\
 \alpha(p(x,y,y),y)=x \;\;\;& 
	\end{align*}
	Then $\VV$ is congruence modular. 
\end{prop}
\proof
First, we check: $\alpha(x,x)=\alpha(p(x,x,x),x)=x$. We  are going now to show that any object $(f,s)$ in a fibre $Pt_Y\VV$ is a canonically endowed with a structure of $\mathrm{Hex}$-algebra (see Definition \ref{VD}), so that this fibre is congruence hyperextensible. By Theorem \ref{teo}, the internal binary operation $\circ$ associated with the term $p$ on any object $(f,s):X\splito Y$ in $Pt_Y\VV$ is defined by $u\circ v=p(u,s(y),v)$ for any $(u,v)\in R[f]$. Then $p(u,sf(u),sf(u))=p(sf(u),sf(u),u)$ shows that it satisfies the first axiom of a $\mathrm{Hex}$-algebra, while the unary  operation defined by $\upsilon(u)=\alpha(u,sf(u))$ satisfies the second one by
$\alpha(p(u,sf(u),sf(u)),sf(u))=u$.
\endproof

The same kind of term-extension obtained from Proposition \ref{exhypo} produces ternary and binary terms satisfying: $\; p(x,x,y)=y \; {\rm and} \; \alpha(p(x,y,y),y)=x$.
But they clearly determine a Mal'tsev term with : $q(x,y,z)=\alpha(p(x,y,z),z)$.

\subsection{Structural outcomes of the characterization theorem}\label{absplit}

Now we shall investigate the structural outcomes of the characterization theorem.  By the point 2) of this theorem and by Proposition \ref{imp2}, in a Gumm category $\EE$, any internal group $(X,\circ,0)$ is necessarily abelian, and is actually determined by a mere subtraction $(X,d,0)$. So, a first immediate outcome is the following one which was already observed in \cite{BGGGum}, but after a long computational way:
\begin{lemma}
	Let $\EE$ be a Gumm category. On any object $X$, there is at most one Mal'tsev operation which is necessarily autonomous.
\end{lemma}
\proof
By Lemma \ref{biter}, a Mal'tsev operation $p$ on the object $X$ in $\EE$ is nothing but a subtraction on the split epimorphism $(p_1^X,s_0^X)$ in the fiber $Pt_X\EE$. Since this fiber is congruence hyperextensible when $\EE$ is a Gumm category, this produces an abelian group structure $(p_1^X,s_0^X)$, which in turm makes $p$ autonomous.
\endproof

But, more than that, the category $\EE$ is immadiately equipped with a \emph{specific class of split epimorphisms}: namely those split epimorphisms $(f,s): X\splito Y$ which are abelian objects in the congruence hyperextensible fiber $Pt_Y\EE$.

So, a split epimorphism $(f,s): X\splito Y$ is \emph{abelian} in $\EE$ if and only if there is a (unique) subtraction on it in the fiber $Pt_Y\EE$, namely if and only if it is endowed with a morphism $d:R[f]\to X$ such that: $d.s_0^f=s.f$ and $d.(1_X,s.f)=1_X$ (i.e. $d(x,x)=sf(x)$ and $d(x,sf(x))=x$); then the morphism $\circ_s: R[f]\to X$ defined by $x\circ_s x'=d(x,d(sf(x),x'))$ necessarily produces an internal group structure with unit $s$ on the  object $f$ of the slice category $\EE/Y$.  We shall denote this class by $APt\EE$ (or $APt$ when  there is no ambiguity about the environment). 
According with what is now usual in the partial Mal'tsev and protomodular
context, it is then natural  to introduce the following:
\begin{defi}
	Let $\EE$ be any Gumm category. An  equivalence relation $R$ on $X$ is said to be \emph{abelian}, when the split epimorphism $(d_0^R,s_0^R)$ is abelian in $Pt_X\EE$, a morphism $f: X\to Y$ is said to be $APt$-special when its kernel equivalence relation  $R[f]$ is abelian. An object $X$ is said to be \emph{affine} when the terminal map $\tau_X$ is $APt$-special or, equivalently, when its indiscrete equivalence relation is $\nabla_X$ abelian.
\end{defi}
Indeed, according to this definition, an affine object $X$ is determined by the existence of a (unique) internal ternary operation $p: X\times  X \times  X \to X$ satisfying the Mal'tsev identities which is then, as we already noticed, necessarily autonomous. So, an equivalence relation $R$ on $X$ is abelian if and only if there is a morphims $p: R[d_0^R]\to X$ in $\EE$ such that:
$$p(aRbRc)Ra; \; p(aRaRb)=b; \; {\rm and}\; p(aRbRb)=a $$
The aim of this  section is to set up the properties of the class $APt\EE$ and its derivations introduced by the previous definition. We get immediately:
\begin{lemma}
The class $APt\EE$ is fibrational, i.e; stable under pullback in $Pt\EE$
\end{lemma}
\proof
Straighforward since, for any map $f: W\to Y$ the base change $f^*:Pt_Y\EE \to Pt_W\EE$ is left exact; so, it preserves any subtraction.
\endproof
Before going any further, we shall need the following observation whose last assertion  makes explicit a very important property which is shared by the (conceptually stronger) Mal'tsev categories:
\begin{prop}
	Let $\EE$ be a Gumm category and the following diagram be any commutative square of split epimorphims in $\EE$ such that $R[f]\cap R[g]=\Delta_W$:
	$$\xymatrix@=3pt{
	W  \ar@<1ex>[rrrr]^{g} \ar@<-1ex>[ddd]_{f}  &&&&  Z  \ar@<-1ex>[ddd]_{\phi_Z} \ar[llll]^{t} \\
		&&&&& && \\
		&&& &&&&&&\\
    X  \ar[rrrr]^{\phi_X}\ar[uuu]_{s} &&&& Y \ar@<1ex>[llll]^{\sigma_X} \ar[uuu]_{\sigma_Z}
	}
	$$
	Then the upward and rightward  square is a pushout in $\EE$. The map $g:W\to Z$, is the unique map $\gamma:W\to Z$ such that $\gamma.s=\sigma_Z.\phi_X$ and $\gamma.t=1_Z$.  This result holds, in particular, when the previous square is a pullback, and shows that the base-change functor $\phi_X^*: Pt_Y\EE\to Pt_X\EE$ along any split epimorphism $\phi_X$ is fully faithful. 
\end{prop}
\proof
It is a consequence of the fact that the fibre $Pt_Y\EE$ is congruence hyperextensible. Suppose given a pair $(h: W\to V,k: Y\to V)$ of morphisms in $\EE$ such that $h.s=k.\phi_X$. Then consider the following morphism of split epimorphisms:
	$$\xymatrix@=3pt{
	W  \ar[rrrr]^{(\phi_X.f,h)} \ar@<-1ex>[ddd]_{f}  &&&&  Y\times V  \ar@<-1ex>[ddd]_{p_Y}  \\
	&&&&& && \\
	&&& &&&&&&\\
	X  \ar[rrrr]^{\phi_X}\ar[uuu]_{s} &&&& Y \ar@{.>}@<1ex>[llll]^{\sigma_X} \ar[uuu]_{(1_Y,k)}
}
$$
This commutative diagram says that the map  $(\phi_X.f,h)$ anihilates the map $s$ in the fiber $Pt_Y\EE$. Since $g$ is the cokernel of $s$ in this fiber, there is a unique factorization $(\phi_Z,\psi): Z\to Y\times V$ such that $(\phi_Z,\psi).g=(\phi_Z,f.h)$ and $(\phi_Z,\psi).\sigma_Z=(1_Y,k)$, namely such that $\psi.g=h$ and $\psi.\sigma_Z=k$. The last sentence is a consequence of Theorem \ref{ffll}.
\endproof
Whence its immediate extension to the regular context: 
\begin{theo}\label{reggumm}
	Let $\EE$ be a regular  Gumm category and  the following diagram be any commutative square of split epimorphims in $\EE$ where 1) $R[f]\cap R[g]=\Delta_W$ and 2) the horizontal arrows are regular epimoprhisms:
	$$\xymatrix@=3pt{
		W  \ar@{->>}[rrrr]^{g} \ar@<-1ex>[ddd]_{f}  &&&&  Z  \ar@<-1ex>[ddd]_{\phi_Z}  \\
		&&&&& && \\
		&&& &&&&&&\\
		X  \ar@{->>}[rrrr]_{\phi_X}\ar[uuu]_{s} &&&& Y  \ar[uuu]_{\sigma_Z}
	}
	$$
	Then the upward and rightward  square is a pushout in $\EE$.  This result holds, in particular, when the previous square is a pullback and $\phi_X$ a regular epimorphism and makes fully faithful any base-change functor $\phi_X^*: Pt_Y\EE\to Pt_X\EE$ when $\phi_X$ is a regular epimorphism.
\end{theo}
\proof
Complete the diagram with the kernel equivalence relations:
	$$\xymatrix@=4pt{
R[g] \ar@<-6pt>[rrrr]_{d_0^g}  \ar@<-1ex>[ddd]_{R(f)} \ar@<6pt>[rrrr]^{d_1^g} 	&&&& W \ar[llll]|{s_0^g} \ar@{->>}[rrrr]^{g} \ar@<-1ex>[ddd]_{f}  &&&&  Z  \ar@<-1ex>[ddd]_{\phi_Z}  \\
	&&&&&&&&& && \\
	&&&&&&& &&&&&&\\
R[\phi_X] \ar@<-6pt>[rrrr]_{d_0^{\phi_X}}   \ar@<6pt>[rrrr]^{d_1^{\phi_X}}\ar[uuu]_{R(s)}	&&&& X \ar[llll]|{s_0^{\phi_X}}  \ar@{->>}[rrrr]_{\phi_X}\ar[uuu]_{s} &&&& Y  \ar[uuu]_{\sigma_Z}
}
$$
When the pair $(f,g)$ is jointly monic, so are the pairs $(R(f),d_0^g)$ and $(R(f),d_1^g)$. Now, according to the previous proposition, any of the upward and rightward left hand side commutative square is a pushout. Consider the two following commutative rectangles which are globally the same:
	$$\xymatrix@=3pt{
	R[\phi_X] \ar@{ >->}[rrrr]^{R(s)}  \ar@{->>}[ddd]_{d_0^{\phi_X}}  	&&&& R[g] \ar@{->>}[rrrr]^{d_1^g} \ar@{->>}[ddd]_{d_0^g}  &&&&  W  \ar@{->>}[ddd]^{g} &&& 		R[\phi_X]\ar@{->>}[rrrr]^{d_1^{\phi_X}} \ar@{->>}[ddd]_{d_0^{\phi_X}}	&&&& X   \ar@{ >->}[rrrr]^{s} \ar@{->>}[ddd]_{\phi_X}&&&& W  \ar@{->>}[ddd]^{g} \\
 &&&& &&&& &&& &&&& &&&&\\
&&&& &&&& &&& &&&& &&&&\\
X \ar@{ >->}[rrrr]_{s}   	&&&& W   \ar@{->>}[rrrr]_{g} &&&& Z  &&& 	X \ar[rrrr]_{\phi_X}	&&&& Y   \ar@{ >->}[rrrr]_{\sigma_Z} &&&& Z
	 }$$
 The two squares on the left hand side are pushouts; so  both rectangles are pushouts. Since the left hand side square of the right hand side rectangle is a pushout, so is the right hand side square of this rectangle.
\endproof

\begin{lemma}
	Let $\EE$ be a Gumm category. Any morphism in $Pt\EE$ between two  split epimorphisms in $APt\EE$ preserves the associated abelian group structures in the fibers. In particular, any morphism between two affine objects preserves the affine structures.
\end{lemma}
\proof
Consider any map in $Pt\EE$ between abelian  objects in the fibers:
$$\xymatrix@=1pt{
	X \ar[dddd]_{f}    \ar[rrrrr]^{\phi}   &&&&&  X'  \ar[dddd]_{f'}  \\
	&&&&& && \\
	&&& &&&&&&\\
	&&& &&&&&&\\
	Y \ar@<-1ex>[uuuu]_{s} \ar[rrrrr]_{\psi} &&&&& Y'  \ar@<-1ex>[uuuu]_{s'}
}
$$
The base change $\psi^*: Pt_{Y'}\to Pt_{Y}$ being left exact produces  an abelian structure on the pullback  $\psi^*(f',s')$ of $(f',s')$ along $\psi$; now, by Proposition \ref{abob}, the induced morphism $(f,s)\to \psi^*(f',s')$ is necessarily a group homomorphism in the congruence hyperextensible fiber $Pt_Y\EE$. As for the last assertion of the lemma, apply the first assertion to the following morphism in $Pt\EE$:
$$\xymatrix@=1pt{
	X\times X \ar[dddd]_{p_1^X}    \ar[rrrrr]^{f\times f}   &&&&&  Z\times Z  \ar[dddd]_{p_1^Z}  \\
	&&&&& && \\
	&&& &&&&&&\\
	&&& &&&&&&\\
	X \ar@<-1ex>[uuuu]_{s_0^X} \ar[rrrrr]_{f} &&&&& Z  \ar@<-1ex>[uuuu]_{s_0^Z}
}
$$
\endproof
From now on, $APt\EE$ can, and will, denote the full subcategory of $Pt\EE$ whose objects are the abelian split epimorphisms. Similarly $AbEqu\EE$ will denote the full subcategory of $Equ\EE$ whose objects are the abelian equivalence relations, $APt\EE/Y$ and $Aff\EE$ will denote the full subcategories of $\EE/Y$ and $\EE$ whose objects are respectively the $APt$-special maps and  the affine objects. 
\begin{coro}\label{reflepi}
	Let $\EE$  be a regular Gumm category. Any base change functor $h^*:Pt_Y\EE \to Pt_{Y'}\EE$ along a regular epimorphism $h$ reflects the abelian split epimorphisms; in other words, given any pullback of split epimorphism:
	$$\xymatrix@=1pt{
		X' \ar[dddd]_{f'}    \ar@{->>}[rrrrr]^{\tilde h}   &&&&&  X  \ar[dddd]_{f}  \\
		&&&&& && \\
		&&& &&&&&&\\
		&&& &&&&&&\\
		Y' \ar@<-1ex>[uuuu]_{s'} \ar@{->>}[rrrrr]_{h} &&&&& Y  \ar@<-1ex>[uuuu]_{s}
	}
	$$
	the split epimorphism $(f,s)$ is abelian as soon as so is $(f',s')$.
	Accordingly, given any fibrant morphism $(f,\tilde f):R\onto S$ of equivalence relations above a regular epimorphism $f:X\onto Y$, the equivalence relation $S$ is abelian if and only if so is $R$.
	More generally, any base change functor $h^*:\EE/Y \to \EE/Y'$ along a regular epimorphism $h$ reflects the $Apt$-special morphisms.
\end{coro}
\proof
Complete the diagram by the horizontal kernel equivalence relations:
$$ \xymatrix@=12pt{
	R[\bar h] \ar@<-1,ex>[dd]_{R(f')} \ar@<-1,ex>[rr]_{p_0^{\bar h}}\ar@<+1,ex>[rr]^{p_1^{\bar h}}
	&& X' \ar@<-1,ex>[dd]_{f'} \ar[ll] \ar[rr]^{\bar h} &&  Y' \ar@<-1,ex>[dd]_{f}\\
	&&&&\\
	R[h] \ar@<-1,ex>[rr]_{p_0^h} \ar@<+1,ex>[rr]^{p_1^h} \ar[uu]_{R(s')} && Y'\ar@{->>}[rr]_{h}\ar[ll]\ar[uu]_{s'}  && Y \ar[uu]_{s} 
}
$$
The right hand side square being a pullback, so are the left hand side ones; accordingly, the split epimorphism $(R(f'),R(s'))$ is abelian, and the two left hand side squares preserve the group structures in the fibers, so that the passage to the quotient gives $(f,s)$ a group structure in the fibre $Pt_Y\EE$.
\endproof

\begin{prop}
Let $\EE$ be a Gumm category. The class $APt\EE$ is point-congruous, namely stable under finite limits in $Pt\EE$ and contains the terminal object $1\splito 1$. 

Accordingly, the subcategories $AbEqu\EE$ and $APt\EE/Y$ (so, in particular, the  subcategory $Aff\EE$) are respectively stable under finite limits in $Equ\EE$ and $\EE/Y$ (resp. in $\EE$). Any category $APt\EE/Y$  (and in particular  $Aff\EE$) is a naturally Mal'tsev category (see Theorem \ref{class1}).
\end{prop}
\proof
The stability under product in $Pt\EE$ is straighforward. Let us check the stability under equalizers. Consider the following diagram where the split epimorphisms $(f,s)$ and $(f',s')$ are abelian and where the two horizontal levels are equalizers:
$$ \xymatrix@=10pt{
I \ar@<-1,ex>[dd]_{\phi} \ar@{ >->}[rr]^{i} && X\ar@<-1,ex>[dd]_{f} \ar@<-1,ex>[rr]_{g_0}\ar@<+1,ex>[rr]^{g_1}
	&& X' \ar@<-1,ex>[dd]_{f'}   \\
	&&&&\\
 J \ar@{ >->}[rr]_{j} \ar[uu]_{\sigma} &&	Y \ar@<-1,ex>[rr]_{h_0} \ar@<+1,ex>[rr]^{h_1} \ar[uu]_{s} && Y' \ar[uu]_{s'} 
}
$$
since the squares indexed by $0$ and $1$ preserve the subtractions in the  fibers, this subtraction naturally extends to the induced split epimorphism $(\phi,\sigma)$. The end of the proposition is a straighforward consequence of the fact that, in any of the categories in question, the finite limits are componentwise.
\endproof

\begin{coro}
Let $\EE$ be a Gumm category. Suppose that the morphisms $g$ and $g.f$ are $APt\EE$-special, then so is $f$.
\end{coro}
\proof
We get the following pullback in $Equ\EE$, where $Y$ is the codomain of $f$:
$$ \xymatrix@=6pt{
	R[f] \ar@{ >->}[dd]_{} \ar@{ >->}[rr]^{} && \Delta_Y \ar@{ >->}[dd]_{} \\
	&&&&\\
	R[g.f] \ar@{ >->}[rr]_{R(f)}  &&  R[g]
}
$$
So $R[f]$ is abelian as soon as so are $R[g.f]$ and $R[g]$.
\endproof

\begin{prop}\label{abab}
	Let $\EE$ be an exact Gumm category and $f: X\onto Y$ any abelian extension, namely any regular $APt$-special morphism $f$. Then the following contruction produces, on the vertical right hand side, an abelian split epimorphism in $\EE$. It is called the \emph{direction} of this extension. This construction applies in particular to the case of  affine objects $X$ with global support, namely such that the terminal map $\tau_X: X \onto 1$ is a regular epimorphism, and, so, it produces an abelian group $\vec A_X$ in $\EE$. 
\end{prop}
\proof
Consider the following diagram where the map $p:R[f]\times_1R[f]\to X$ is defined by $p(x,y,z)=\partial(x,\partial(y,z))$. Since, in a Gumm context, this ternary operation is autonomous, it makes pullback any commutative square on the left hand side:
$$ \xymatrix@=15pt{
	R[f]\times_1 R[f] \ar@<-1,ex>[dd]_{p_0^X}\ar@<+1,ex>[dd]^{(p,p_1^X.p_1^X)} \ar@<-1,ex>[rr]_{(p_0^X.p_1^X,p)}\ar@<+1,ex>[rr]^{p_2^X}
	&& R[f] \ar@<-1,ex>[dd]_{p_0^f}\ar@<+1,ex>[dd]^{p_1^f} \ar[ll] \ar@{.>>}[rr]^{q_f} && \vec A_f \ar@{.>>}@<+1,ex>[dd]^{\psi_f}\\
	&&&&\\
	R[f] \ar@<-1,ex>[rr]_{p_0^f} \ar@<+1,ex>[rr]^{p_1^f} \ar[uu]_{} && X\ar@{->>}[rr]_{f}\ar[ll]\ar[uu]_{}  && Y \ar@{.>}[uu]^{0_f} 
}
$$
Let $q_f$ be the quotient map of the upper horizontal equivalence relation. Since the left hand  side commutative squares are pullbacks, so is the right hand side one in the exact category $\EE$ thanks to the Barr-Kock Theorem. Then, as in any category $\EE$, the autonomous law  $p$ gives $(0_f,\psi_f)$ an abelian  structure defined by $\vec{ab}+\vec{bc}=\vec{ac}$, see \cite{B3}.
\endproof

\noindent\textbf{Baer sum of abelian extensions with a given direction}

So, again according to \cite{B3} and similarly to the exact Mal'tsev context, there is, when $\EE$ is an exact Gumm category, an abelian group structure on the set $Ext_{\vec A}$ of isomorphic classes of extensions $f:X\onto Y$ with abelian kernel relation $R[f]$ having as direction a given abelian group $(\psi,0): \vec A\splito Y$ in $Pt_Y\EE$, the binary operation of this group being produced by the so-called \emph{Baer sum}.

\subsection{Preconnectors, internal groupoids and categories}\label{cat1}

In this section we shall be interested by the nature of the internal groupoids and categories in the Gumm context. Some results of this section were already stated in \cite{BGGGum} where they were checked in a computational and combinatorial way quite different than the structural one followed in this section and resulting from the first proposition. The Propositions \ref{new} and \ref{antepen} introduce  important structural precisions about the nature of these internal categories and groupoids. 

Consider any preorder $T$ on an object $Y$ in the category $\EE$:
$$\xymatrix{  T\times_YT \ar@<-6pt>[r]_{d_0^T}\ar[r]|<<<{d_1^T}  \ar@<6pt>[r]^{d_2^T} 	& T \ar@<-6pt>[r]_{d_0^T} \ar@<6pt>[r]^{d_1^T}& Y \ar[l]|{s_0^T} 
}
$$
Let $\delta_0: T\times_1 X \to X$ be a map making commutative the right hand side square  indexed by $0$ in the  following diagram where the  square indexed  by $1$ is a pullback:
$$\xymatrix@=2pt{
	T\times_YT\times_1X \ar[dddd]_{h_2} \ar@<-8pt>@{.>}[rrrr]_{T\times_1\delta_0} \ar[rrrr]|{\delta_1^h}  \ar@<8pt>[rrrr]^{\delta_2^h}  &&&&  T\times_1X \ar@<-6pt>@{.>}[rrrr]_{\delta_0}  \ar[dddd]_{h_1} \ar@<6pt>[rrrr]^{\delta_1^h}  &&&&  X  \ar[dddd]^{h} \ar[llll]|{\sigma_0^h} \\
	&&&&& && \\
	&&& &&&&&&\\
	&&& &&&&&&\\
	T\times_YT \ar@<-8pt>[rrrr]_{d_0^T} \ar[rrrr]|{d_1^T}  \ar@<8pt>[rrrr]^{d_2^T} 	&&&& T \ar@<-8pt>[rrrr]_{d_0^T} \ar@<8pt>[rrrr]^{d_1^T} &&&& Y \ar[llll]|{s_0^T} 
}
$$ 
Then the pair $(\delta_0,\delta_1^h)$ determines a relation on $X$ since $T$ is a relation on $Y$; and $h$ becomes a vertical morphism of relations. If we suppose:  (1) $\delta_0.\sigma_0^h=1_X$, the upper horizontal relation becomes reflexive.

As  usual, a \emph{left action} of the preorder $T$ on the map $h: Z\to Y$ (or on the object $h$ in the slice category $\EE/Y$) is given by a map $\delta_0: T\times_1 Z \to Z$ making commutative the square  indexed by $0$
which, moreover, satisfies (1)  and (2): $\delta_0.\delta_1^h=\delta_0.T\times_1\delta_0$. Then the higher horizontal part on this diagram produces a preorder on $X$ and the morphism $h$ a fibrant morphism of preorders. A \emph{right action} of $T$ on $h$ is a left action of the dual preorder $T^{op}$.
 \begin{prop}\label{action}
 	Given any Gumm category $\EE$, a morphism $\delta_0$ making commutative the above square indexed by $0$ determines a left action of the preorder $T$ on $h$ if and only if it satisfies axiom (1).
 	
 	Under this condition, when  $T$ is an equivalence relation, the upper horizontal diagram is an equivalence relation as well.  When $h$ is split by $t$ and $t_1$ denotes the induced splitting of $h_1$, then there is at most one action of the equivalence relation $T$ on $h$ which preserves the splitting, namely such that $\delta_0.t_1=d_0.t$.
 \end{prop}
 \proof
 Using the Yoneda embedding, it is enough to check it in $Set$. Consider the following situation in the set $T\times_1X$:
 $$\xymatrix@=7pt{
 	(bTc,x) \ar[rr]^{R[\delta_1^h]} \ar[dd]^{R[d_0^T.h_1]}  \ar@(l,l)[dd]_{R[\delta_0]} && (aTc,x)\ar[dd]_{R[d_0^T.h_1]} \ar@(r,r)@{.>}[dd]^{R[\delta_0]}\\
 	&&&\\
 	(bTb,\delta_0(bTc,x)) \ar[rr]_{R[\delta_1^h]} && (aTb,\delta_0(bTc,x))
 }
 $$
 We have $R[h_1]\cap R[\delta_1^h]=\Delta_{T\times_1X}$ since the square indexed by $1$ is a pullback. Since $T$ is a relation, we still get $R[d_0^T.h_1]\cap R[\delta_1^h]=\Delta_{T\times_1X}$. Accordingly  $R[d_0^T.h_1]\cap R[\delta_1^h]\subset R[\delta_0]$. So, by the Shifting Property, we get the doted vertical arrow which means $\delta_0((aTc,x))=\delta_0(aTb,\delta_0(bTc,x))$ which is (2).
 
 That the upper preorder is an equivalence relation when so is $T$ is straightforward. Let us check now that when $T$ is an equivalence relation, and $\delta_0$ is an  action, then $R[h_1]\cap R[\delta_0]=\Delta_{T\times_1X}$, i.e.  the pair $(h_1,\delta_0)$ is jointly monic. Suppose $\delta_0(bTc,x)=\delta_0(bTc,x')$; then $x=\delta_0(cTb,\delta_0(bTc,x))=\delta_0(cTb,\delta_0(bTc,x'))=x'$. Now, saying that the action preserves the splitting $t$ of $h$ means that $\delta_0(bTc,t(c))$ $=t(b)$; i.e. it means  that the following commutative diagram is a punctual span in the  fibre $Pt_Y\EE$, actually a punctual relation according to our first observation about $(h_1,\delta_0)$:
 $$\xymatrix@=1pt{
 	T\times_1X \ar@<-6pt>[rrrrrr]_{\delta_0}  \ar@<-6pt>[dddd]_{h_1}  &&&&&&  X  \ar\ar@<-4pt>[dddd]_{h} \ar[llllll]_{\sigma_0^h} \\
 	&&&&&&&& \\
 	&&&&&&&&\\
 	&&&&&&&&\\
 	T \ar@<-6pt>[rrrrrr]_{d_0^T} \ar[uuuu]_{t_1} &&&&&& Y \ar[llllll]_{s_0^T} \ar[uuuu]_{t} 
 }
 $$
 Consequently the morphism $\delta_0$ is the unique one such that $\delta_0.\sigma_0^h=1_X$ and $\delta_0.t_1=t.d_0^T$. 
 \endproof
Given any pair $(R,S)$ of equivalence relations on an object $X$, a \emph{preconnector} between them (see Section \ref{conn}) is a subternary operation $p:R\times_1S\to X$ such that $xSp(xRySz)Rz$:
$$\xymatrix@=0,5pt{
	x  \ar[dddd]_{R} \ar@{.>}[rrrr]^{S}  &&&&  p(xRySz)  \ar@{.>}[dddd]^{R}  \\
	&&&&& \\
	&&&&&\\
	&&&&&\\
	y \ar[rrrr]_{S} &&&& z 
}
$$
It produces a splitting $\pi_p$ of $\zeta_{(R,S)}: R\square S \to R\times_1S$ which, in the regular context, asserts the permutation of this pair, see Proposition \ref{keyperm} below.
\begin{lemma}\label{preconn}
Let $\EE$ be a Gumm category. Suppose that, in addition, the preconnector $p$ on $(R,S)$ satisfies $p(xRySy)=x$,  then we get: $p(xRySz)=p(x'RySz)$ $\iff x=x'$. So, we have: $R[d_1^S]\cap R[p]=\Delta_{R\times_1S}$.
\end{lemma}
\proof
Consider  the following situation:
 $$\xymatrix@=5pt{
	xRySz \ar[rr]^{R[\bar d_0^R]} \ar[dd]^{R[\bar d_1^S]} \ar@(l,l)[dd]_{R[p]} && xRySy\ar[dd]_{R[\bar d_1^S]} \ar@{.>}@(r,r)[dd]^{R[p]}\\
	&&&\\
	x'RySz \ar[rr]_{R[\bar d_0^R]} && x'RySy
}
$$
Since $R[\bar d_0^R]\cap R[\bar d_1^S]=\Delta_{R\times_1S}$, we get the right hand side bended arrow which, according to our asumption, implies $x=x'$. 
\endproof

\begin{theo}[\cite{BGGGum}]\label{affmal}
	Given any Gumm category $\EE$ and any pair $(R,S)$ of equivalence relations on an object $X$, there is, between them, at most one preconnector satisfying the  Mal'sev identities: $(i): p(aRaSb)=b$ and $(ii): p(aRcSc)=c$. These two axioms imply the missing one which completes the definition of a \emph{connector}, see Section \ref{conn}.  When $f: X\to Y$ is a morphism of pairs of equivalence relations $(R,S)\to (R',S')$ which are respectively preconnected with $(i)$ and $(ii)$, then the map $f$ preserves the preconnectors.
\end{theo}
\proof 
Consider the following left hand side double equivalence relation associated with the square construction, see Section \ref{square}:
$$\xymatrix@=4pt{
	R\square S \ar@<-2ex>[dddd]_{\delta_0^R} \ar@<-2ex>[rrrrr]_{\delta_0^S}    \ar@<2ex>[dddd]^{\delta_1^R} \ar@<2ex>[rrrrr]^{\delta_1^S}  &&&&&  S  \ar@<-2ex>[dddd]_{d_0^S} \ar[lllll]|{\sigma_0^S}\ar@<2ex>[dddd]^{d_1^S} &&&&& 	R\times_1S  \ar@<-2ex>[dddd]_{\bar d_0^R} \ar@<-2ex>[rrrrr]_{\delta_0^S.\pi_p}    \ar@<2ex>@{.>}[dddd]^{\delta_1^R.\pi_p} \ar@<2ex>[rrrrr]^{\bar d_1^S}  &&&&&  S  \ar@<-2ex>[dddd]_{d_0^S} \ar[lllll]|{\bar s_0^S}\ar@<2ex>@{.>}[dddd]^{d_1^S}\\
&&&&& && &&&&&\\
	&&& &&&&&&&&&&&\\
	&&& &&&&&&&&&&&\\
	R \ar[uuuu]|{\sigma_0^R} \ar@<-2ex>[rrrrr]_{d_0^R} \ar@<2ex>[rrrrr]^{d_1^R}  &&&&& X \ar[lllll]|{s_0^R} \ar[uuuu]|{s_0^S} &&&&& 	R \ar[uuuu]|{\bar s_0^R} \ar@<-2ex>[rrrrr]_{d_0^R} \ar@<2ex>[rrrrr]^{d_1^R}  &&&&& X \ar[lllll]|{s_0^R} \ar[uuuu]|{s_0^S}
}
$$
With respect to the splitting $\pi_p$ of $\zeta_{R,S}: R\square S \to R\times_1S$ described above, 
the previous identities become on the right hand  side: $(i): \delta_1^R.\pi_p.\bar s_0^R=1_R$ or, equivalently $\delta_0^S.\pi_p.\bar s_0^S=s_0^R.d_1^S$; and $(ii): \delta_0^S.\pi_p.\bar s_0^S=1_S$ or, equivalently, $\delta_0^S.\pi_p.\bar s_0^R=s_0^S.d_0^R$. So, now, concerning the right hand side diagram, by $(ii)$, we get the upper horizontal reflexive relation and by $(i)$ the plain vertical split epimorphism between the upper horizontal reflexive relation and the lower horizontal equivalence relation. This vertical split epimorphism fullfils the conditions of the previous proposition; so, $\delta_0^S.\pi$ is the unique map
 $\delta: R\times_1S\to S$ such that ($\delta.\bar s_0^S=1_S$ and) $\delta.\bar s_0^R=s_0^S.d_0^R$. On the other hand, since $\delta_1^S.\pi$ is necessarily $d_1^S$ and the pair $(\delta_0^S,\delta_1^S)$ is jointly monomorphic, the splitting $\pi$ is unique as well which implies the uniqueness of the preconnector.
 
 According to the Proposition \ref{action}, since we are in a Gumm category, the axiom $(ii)$ is enough to produce a  left action of the equivalence relation $R$ on the split epimorphim $(d_0^S,s_0^S)$. So, we get for free the second axiom of a left action, namely $p(aRbSp(bRcSd))=p(aRcSd))$, i.e. the half of Axiom 3 for a connector. The second half is obtained by inverting the role of $R$ and $S$. 
 
 Let $f: X\to Y$  be a morphism of pairs of equivalence relations $(R,S)\to (R',S')$ which are respectively preconnected with conditions $(i)$ and $(ii)$, and $q$ and $q'$ denote the respective ternary operation associated with these preconnectors.
 Consider now the following situation $R\times_1S$:
 $$\xymatrix@=7pt{
 	aRbSc \ar[rr]^{R[d_1^S]} \ar[dd]^{R[q]}  \ar@(l,l)@{.>}[dd]_{R[\psi]} && bRbSc\ar[dd]_{R[q]} \ar@(r,r)[dd]^{R[\psi]}\\
 	&&&\\
 	q(aRbSc)RcSc \ar[rr]_{R[d_1^S]} && cRcSc
 }
 $$
 where $\psi(aRbSc)=q'(f(a)R'f(b)S'f(c))$. Since, by the previous lemma, we have $R[d_1^S]\cap R[q]=\Delta_{R\times_1S}\subset R[\psi]$, we get the bended left hand side arrow which means: $q'(f(a)R'f(b)S'f(c))=fq(aRbSc)$.
\endproof

When such a connector exists for a pair $(R,S)$, we shall say, as in the Mal'tsev context, that the pair of equivalence relations centralizes each other. From Theorem \ref{affmal} that we get immediately:
\begin{coro}[\cite{BGGGum}]\label{intgrds}
	Let $\EE$ be a Gumm category. On any reflexive graph $X_{\bullet}$:
	$$\xymatrix@=2pt{
	  X_1 \ar@<-6pt>[rrrr]_{d_0^{X_{\bullet}}}   \ar@<6pt>[rrrr]^{d_1^{X_{\bullet}}}  &&&&  X_0   \ar[llll]|{s_0^{X_{\bullet}}}  
	}
	$$ 
	 there is at most one groupoid structure. The forgetful functor $Grd\EE \to RGh\EE$ from the internal groupoids to the reflexive graphs in $\EE$ is full (and faithful, as in any case).
\end{coro}
\proof
Straighforward, since a graph is a groupoid if and only if the pair of equivalence relations  $(R[d_0^{X_{\bullet}}],R[d_1^{X_{\bullet}}])$ is connected, see \cite{JP} and Section \ref{conn} below. 
\endproof 

Let us denote $(\;)_0: Grd\EE \to \EE$ the forgetful functor associating with any groupoid $X_{\bullet}$ its "object of objects" $X_0$. We can the add the following important precisions:
\begin{prop}\label{antepen}
	Let $\EE$ be any Gumm category. Then any groupoid $X_{\bullet}$ is abelian. Accordingly, any fiber $Grd_X\EE$ of the fibration $(\;)_0$ is a naturally Mal'tsev category. Furthermore, this fiber is such that any base-change functor with respect to its fibration of point is fully faithful.
\end{prop}
\proof
Given any category $\EE$ and any internal groupoid $X_{\bullet}$, the map $(d_0^{X_{\bullet}},d_1^{X_{\bullet}}): X_1\to X_0\times X_0$ is canonically endowed with an associative Mal'tsev operation in the slice category $\EE/(X_0\times X_0)$ given by the map $\chi:R_2[(d_0^{X_{\bullet}},d_1^{X_{\bullet}})]\to  X_1$ defined by $\chi(\alpha,\beta,\gamma)=(\alpha.\beta^{-1}.\gamma)$. A groupoid is defined as \emph{abelian} when this associative Mal'tsev operation is autonomous. This is necessarily the case by Theorem \ref{affmal} and since, according to the previous corollary, the fiber $Grd_X\EE$ is a full subcategory of the fiber $RGh_X\EE$ of reflexive graphs with object of objects $X$, which is itself a Gumm category as the coslice category $s_0^X/(\EE/(X\times X))$:
$$\xymatrix@=7pt{
	& G   \ar[ddr]^{(d_0,d_1)}\\
	&&&\\
	X \ar@{ >->}[rr]_{s_0^{X}} \ar@{ >->}[ruu]^{s_0^{G}} && X\times X
}
$$
So, since any of its objects is affine, the fiber $Grd_X\EE$ is a naturally Mal'tsev category. The last point is  shown in \cite{Bab} where is introduced a classification table of the different possible notions of non-pointed "additive categories" and whose one level is the property in question.
\endproof

Let us now move on to the internal categories:

\begin{prop}[\cite{BGGGum}]\label{intcat}
	Let $\EE$ be a Gumm category. On any reflexive graph $X_{\bullet}$ 
	there is at most one internal category structure which is necessarily right and left cancellable. For that it is enough to get the composition map $d_1:X_1\times_0X_1 \to X_1$ (where $d_1(\alpha,\beta)$ denotes $\beta.\alpha$ in order to respect the simplicial notations) such that, for any pair $(\alpha,\beta)$ of composable arrows in $X_1$, we get: 
	\begin{align*}
 1)\; d_1(s_0d_0(\beta),\beta)=\beta\;\;\;\;   & \;\; 1')\; d_1(\alpha,s_0d_1(\alpha))=\alpha\\
 2)\;  d_0^{X_{\bullet}}(d_1(\alpha,\beta))=d_0^{X_{\bullet}}(\alpha)\;  & \;\; 2')\; d_1^{X_{\bullet}}(d_1(\alpha,\beta))=d_1^{X_{\bullet}}(\beta)
  \end{align*}
  In other words, we get the associativity for free.
 Any morphism between the underlying graphs of two internal categories in necessarily a functor. Any fiber $Cat_{X_0}\EE$ is a Gumm category.
\end{prop}
\proof
Consider the following diagram, where the rightward and downward left hand side square of plain arrows is underlying a pullback of split epimorphisms: 
$$\xymatrix@=2pt{
  X\times_0X \ar@<-8pt>@{.>}[rrrrrr]_{d_1}  \ar@<-6pt>[dddd]_{d_0^{X_{\bullet}}} \ar@<8pt>[rrrrrr]^{d_2^{X_{\bullet}}}  &&&&&&  X_1  \ar@<-6pt>[dddd]_{d_0^{X_{\bullet}}} \ar[llllll]|{s_1^{X_{\bullet}}} \ar[rrrrrrr]^{d_1^{X_{\bullet}}} &&&&&&& X_0\\
	&&&&&& \\
    &&&&&&\\
	&&&&&&\\
  X_1 \ar@<-8pt>@{.>}[rrrrrr]_{d_0^{X_{\bullet}}} \ar@<8pt>[rrrrrr]^{d_1^{X_{\bullet}}}\ar[uuuu]_{s_0^{X_{\bullet}}} &&&&&& X_0 \ar[llllll]|{s_0^{X_{\bullet}}} \ar[uuuu]_{s_0^{X_{\bullet}}}
}
$$
and the map $d_1$ depicts the composition map when  it exists. Axiom 2) is equivalent to the commutation of the vertical square with dotted arrows: $d_0^{X_{\bullet}}.d_1=d_0^{X_{\bullet}}.d_0^{X_{\bullet}}$, while Axiom 1) is equivalent to $d_1.s_1^{X_{\bullet}}=1_{X_1}$; so, the axioms 1) and 2) delineate a kind of ''action" of the graph $X_{\bullet}$ on the map $d_0^{X_{\bullet}}: X_1\to X_0$. Axiom 1') is equivalent to $d_1.s_0^{X_{\bullet}}=1_{X_1}$ while Axiom 2') says that the  map $d_1^{X_{\bullet}}$ coequalizes the pair $(d_1,d_2^{X_{\bullet}})$. 

By the Shifting Lemma, Axiom 1) implies that the pair $(d_0^{X_{\bullet}},d_1)$ is jointly monic so  that the composition of arrows is right cancellable: considering the following situation in $X_1\times_0X_1$ depicting $\beta.\alpha=\beta'.\alpha$:
$$\xymatrix@=7pt{
	(\alpha,\beta) \ar[rr]^{R[d_2^{X_{\bullet}}]} \ar[dd]^{R[d_0^{X_{\bullet}}]}  \ar@(l,l)[dd]_{R[d_1]} && (s_0d_1(\alpha),\beta) \ar[dd]_{R[d_0^{X_{\bullet}}]} \ar@(r,r)@{.>}[dd]^{R[d_1]}\\
	&&&\\
	(\alpha,\beta') \ar[rr]_{R[d_2^{X_{\bullet}}]} && (s_0d_1(\alpha),\beta')
}
$$
the right hand side induced bended arrows means $\beta=\beta'$. By symmetry Axiom 1') implies that the pair $(d_2^{X_{\bullet}},d_1)$ is jointly monic; so, 1) the composition of arrows is left cancellable as well, and 2) the  upper horizontal graph is actually a reflexive relation.

From that we first get the uniqueness of a composition map $d_1$ satisfying 1), 1') and 2): considering the following situation in $X_1\times_0X_1$:
$$\xymatrix@=7pt{
	(\alpha,\beta) \ar[rr]^{R[d_0^{X_{\bullet}}]} \ar[dd]^{R[d_1]}  \ar@(l,l)@{.>}[dd]_{R[d'_1]} && (\alpha,s_0d_1(\alpha))) \ar[dd]_{R[d_1]} \ar@(r,r)[dd]^{R[d'_1]}\\
	&&&\\
	(s_0d_0(\alpha),\beta.\alpha) \ar[rr]_{R[d_0^{X_{\bullet}}]} && (s_0d_0(\alpha),\alpha)
}
$$
where $(s_0d_0(\alpha),\beta.\alpha)$ is a composable pair by axiom 2) and where the  vertical square with $R[d_1]$ is a consequence of 1) and 1'), we get $d'_1=d_1$ (namely the left hand side bended arrow) provided that we have the right hand side bended arrow, namely when $d'_1$ satisfies 1) and 1') as well.

Finally we shall get the associativity of the composition in the following way: first, let us extend the previous diagram by the kernel equivalence relations of the maps $d_1^{X_{\bullet}}$ and $d_2^{X_{\bullet}}$, where the indexation of the left hand part of the diagram  respects the simplicial ones:
$$\xymatrix@=2pt{
R[d_2^{X_{\bullet}}] \ar@<-8pt>[rrrrrr]_{d_2}   \ar@<8pt>[rrrrrr]^{d_3} \ar@<-6pt>[dddd]_{R(d_0^{X_{\bullet}})} &&&&&&	X\times_0X  \ar[llllll]|{s_2} \ar@<-12pt>@{.>}[rrrrrr]_{d_1}  \ar@<-6pt>[dddd]_{d_0^{X_{\bullet}}} \ar@<4pt>[rrrrrr]^{d_2^{X_{\bullet}}}  &&&&&&  X_1  \ar@<-2pt>[dddd]^{d_0^{X_{\bullet}}} \ar@<4pt>[llllll]|{s_1^{X_{\bullet}}} \\
 &&&&&&	&&&&&& \\
 &&&&&&	&&&&&&\\
 &&&&&&	&&&&&&\\
R[d_1^{X_{\bullet}}] \ar@<-8pt>[rrrrrr]_{d_1}   \ar@<8pt>[rrrrrr]^{d_2} &&&&&&	X_1 \ar[llllll]|{s_1} \ar@<-12pt>@{.>}[rrrrrr]_{d_0^{X_{\bullet}}} \ar@<4pt>[rrrrrr]^{d_1^{X_{\bullet}}}  &&&&&& X_0 \ar@<4pt>[llllll]|{s_0^{X_{\bullet}}} 
}
$$
From $(d_0^{X_{\bullet}},d_1)$ jointly monic, we get $(R(d_0),d_1.d_2)$ jointly monic; 
now, starting with a triple $(\alpha,\beta,\gamma)$ of composable arrows in $X_1$, let us consider the following situation in $R[d_2^{X_{\bullet}}]$ where we need Axiom 2') to make composable the pair $(\beta.\alpha,\gamma)$:
$$\xymatrix@=7pt{
	(\beta,\beta.\alpha,\gamma) \ar[rr]^{R[R(d_0)]} \ar[dd]^{R[d_1.d_2]}  \ar@(l,l)@{.>}[dd]_{R[d_1.d_3]} && (\beta,\beta.\alpha,d_0(\gamma)) \ar[dd]_{R[d_1.d_2]} \ar@(r,r)[dd]^{R[d_1.d_3]}\\
	&&&\\
	(s_0d_0(\beta),\alpha,\gamma.\beta) \ar[rr]_{R[R(d_0)]} && (s_0d_0(\beta),\alpha,\beta)
}
$$
where the left hand side dotted bended arrow means  $\gamma.(\beta.\alpha)=(\gamma.\beta).\alpha$. Actually, this associativity turns the following reflexive relation into a transitive one:
$$\xymatrix@=2pt{
	X\times_0X \ar@<-8pt>[rrrrrr]_{d_1}   \ar@<8pt>[rrrrrr]^{d_2^{X_{\bullet}}} &&&&&& X_0 \ar[llllll]|{s_1^{X_{\bullet}}}
}$$

Now, given a morphism of reflexive graphs:
$$\xymatrix@=2pt{
	X_1 \ar[rrrrrr]^{f_1}  \ar@<-6pt>[dddd]_{d_0^{X_{\bullet}}} \ar@<6pt>[dddd]^{d_1^{X_{\bullet}}}  &&&&&&  Y_1 \ar@<6pt>[dddd]^{d_1^{Y_{\bullet}}}  \ar@<-6pt>[dddd]_{d_0^{Y_{\bullet}}}  \\
	&&&&&& \\
	&&&&&&\\
	&&&&&&\\
	X_0 \ar[uuuu]|{s_0}\ar[rrrrrr]_{f_0}  &&&&&& Y_0 \ar[uuuu]|{s_0} 	
}
$$
consider the following situation in $X_1\times_0X_1$:
$$\xymatrix@=7pt{
	(\alpha,\beta) \ar[rr]^{R[d_0^{X_{\bullet}}]} \ar[dd]^{R[d_1]}  \ar@(l,l)@{.>}[dd]_{R[d_1,(f_1\times_0f_1)]} && (\alpha,s_0d_1(\alpha)) \ar[dd]_{R[d_1]} \ar@(r,r)[dd]^{R[d_1,(f_1\times_0f_1)]}\\
	&&&\\
	(s_0d_0(\alpha),\beta.\alpha) \ar[rr]_{R[d_0^{X_{\bullet}}]} && (s_0d_0(\alpha),\alpha)
}
$$
where the left hand side bended arrow means: $f_1(\beta).f_1(\alpha)=f_1(\beta.\alpha)$. For the last point, see the end of the proof of the last Proposition.
\endproof
The following precisions about internal categories are new and important:
\begin{prop}\label{new}
	Let $\EE$ be a Gumm category. On any pointed object $(X,e)$ in $\EE$ there is at most one unitary magma with unit $e$ which is necessarily associative, left cancellable and commutative. Any pointed morphism $f: (X,e)\to (X',e')$ between such unitary magmas is a monoid homomorphism. Any internal category $X_{\bullet}$ 
 has necessarily the monoids of its endomorphisms left cancellable and commutative. 
\end{prop}
\proof
For the first assertion, apply the previous proposition to the following reflexive graph:
$$\xymatrix@=2pt{
	X \ar@<-6pt>[rrrr]_{\tau_X}   \ar@<6pt>[rrrr]^{\tau_X}  &&&&  1   \ar[llll]|{e}  
}
$$
Denote by $m=d_1$ the composition with unit $e$. The uniqueness of the internal category structure implies in this case that $m=m^{op}$ which means that the monoid is  commutative.

More generally, given any internal category $X_{\bullet}$, let us define $EndX_{\bullet}$ by the following pullback in the category of internal categories $Cat\EE$:
$$\xymatrix@=7pt{
	EndX_{\bullet} \ar@{ >->}[rr]^{} \ar[dd]^{}   && X_{\bullet} \ar[dd]_{} \\
	&&&\\
	\Delta_{X_0} \ar@{ >->}[rr] \ar@{ >.>}[rruu]&& \nabla_{X_0}
}
$$
The left hand side arrow is split because of the diagonal inclusion; actually this split epimorphism in $Cat\EE$ determines a monoid in the Gumm slice category $\EE/X_0$; it is then left cancellable and commutative.
\endproof
So, the internal categories in a Gumm categories are quite specific: they are left and right cancellable and such that each monoid of endomorphims is commutative and left cancellable: it is what it is remaining, on internal categories, of Proposition \ref{antepen} for internal groupoids.

It is clear that when $\EE$ is a Gumm category, any functor category $\mathcal{F}(\CC,\EE)$ is a Gumm category as well since the limits in the functor categories are pointwise. Accordingly, when $\EE$ is a Gumm category so is the category $RGh\EE$ of internal reflexive graphs. By Corollary \ref{intgrds} and Proposition \ref{intcat}, the forgetful functors:
$$Grd\EE \to Cat\EE \to RGh\EE$$
are full, faithful and left exact, so that the two first categories appear as full and faithful subcategories of the last one in which they are stable under finite limits. Whence immediately:
\begin{prop}
	When $\EE$ is a Gumm category, so are $Grd\EE$ and $Cat\EE$.
\end{prop}

 \subsection{Discrete fibrations and cofibrations}\label{cat2}
 
A discrete fibration (resp. cofibration) is a functor between internal categories:
$$\xymatrix@=14pt{
	X_2  \ar@<2ex>[rr]^{d_{2}^{X_{\bullet}}}\ar[dd]_{h_2} \ar[rr]|{d_{1}^{X_{\bullet}}} \ar@<-2ex>[rr]_{d_{0}^{X_{\bullet}}} && X_1 \ar@<2ex>[rr]^{d_{1}^{X_{\bullet}}} \ar@<-2ex>[rr]_{d_{0}^{X_{\bullet}}}  \ar[dd]_{h_1}&&
	X_0   \ar[ll]|{s_0^{X_{\bullet}}} \ar[dd]^{h_0}\\
	&&&& \\
	Y_2  \ar@<2ex>[rr]^{d_{2}^{Y_{\bullet}}} \ar[rr]|{d_{1}^{Y_{\bullet}}} \ar@<-2ex>[rr]_{d_{0}^{Y_{\bullet}}} && Y_1 \ar@<2ex>[rr]^{d_{1}^{Y_{\bullet}}} \ar@<-2ex>[rr]_{d_{0}^{Y_{\bullet}}} &&
	Y_0   \ar[ll]|{s_0^{Y_{\bullet}}}
	  }
$$
for which the right hand side diagram indexed by $1$ (resp. by $0$) is a pullback.
By \cite{J} for instance, we know that the category $Difib_{Y_{\bullet}}$ of discrete fibration above $Y_{\bullet}$ is the category of algebras of a left exact monad $(T,\lambda,\mu)$ on the slice category  $\EE/Y_0$ determined by the internal category $Y_{\bullet}$. Whence:
\begin{prop}
When $\EE$ is Gumm category, so is any category $Difib_{Y_{\bullet}}$. By duallity so is the category $Dicofib_{Y_{\bullet}}$ of discrete cofibration above $Y_{\bullet}$.
\end{prop}
We are going now to make explicit that, in the  Gumm context, the objects of $Difib_{Y_{\bullet}}$ are very simple ones:
\begin{prop}
	Let $\EE$ be a Gumm category and $Y_{\bullet}$ an internal category in $\EE$. Consider a morphism of reflexive graphs where the vertical commutative square indexed by $1$ is pullback:
	$$\xymatrix@=1,5pt{
		X_1 \ar@<-8pt>[rrrrrr]_{d_1^{X_{\bullet}}}  \ar[ddddd]_{h_1} \ar@<8pt>[rrrrrr]^{d_0^{X_{\bullet}}}  &&&&&&  X_0  \ar[ddddd]^{h_0} \ar[llllll]|{s_0^{X_{\bullet}}}  \\
		&&&&&& \\
		&&&&&&\\
		&&&&&&\\
		&&&&&&\\
		Y_1 \ar@<-8pt>[rrrrrr]_{d_0^{Y_{\bullet}}} \ar@<8pt>[rrrrrr]^{d_1^{Y_{\bullet}}} &&&&&& Y_0 \ar[llllll]|{s_0^{Y_{\bullet}}} 
	}
	$$
	Then the graph $X_{\bullet}$ is actually underlying an internal category and, accordingly, the functor $h_{\bullet}$ a discrete fibration. This functor, in addition, is such that, given any arrow $\beta \in Y_1$, there is, above it, at most one arrow $\alpha \in X_1$ with a given domain.
	\end{prop}
\proof
Since the square indexed by $1$ is a pullback, any composable pair $(\alpha,\beta)\in X_1$  produces a unique map $\bar d_1(\alpha,\beta)$ above $d_1(h_1(\alpha),h_1(\beta))$ with codomain $d_1(\beta)$.
It is straightforward to check that $\bar d_1$ satifies the axioms  1),1') and 2') of Proposition \ref{intcat} since $d_1$ satisfies them. It remains to check 2). For that
consider the following situation in $X_1\times_0X_1$, where we get the right hand side bended arrow  by 1) and 1'):
$$\xymatrix@=6pt{
	(\alpha,\beta) \ar[rrrr]^{R[d_2^{X_{\bullet}}]} \ar[ddd]^<<<<{R[d_1.(h_1\times h_1)]}  \ar@(l,l)@{.>}[ddd]_{R[d_0^{X_{\bullet}}.d_0^{X_{\bullet}}]} &&&& (s_0d_0(\beta),\beta) \ar[ddd]_>>>>{R[d_1.(h_1\times h_1)]} \ar@(r,r)[ddd]^{R[d_0^{X_{\bullet}}.d_0^{X_{\bullet}}]}\\
	&&&&\\
	&&&&\\
	(\bar d_1(\alpha,\beta),s_0d_1(\beta)) \ar[rrrr]_{R[d_2^{X_{\bullet}}]} &&&& (\beta,s_0d_1(\beta))
}
$$
The pair $(d_1.(h_1\times h_1),d_2^{X_{\bullet}})$ is jointly monic, since so is $(d_1,d_2^{Y_{\bullet}})$, the composition being left cancellable  in the category $Y_{\bullet}$. Accordingly, we get the right hand side bended arrow which means $d_0\bar d_1(\alpha,\beta)=d_0(\alpha)$, namely 2). So, the graph $X_{\bullet}$ is underlying a category, $h_{\bullet}$ becomes a functor which is a discrete fibration since our starting assumption is that the square indexed by $1$ is a pullback. The last assertion about this functor  means that the pair $(h_1,d_0^{X_{\bullet}})$ is jointly monic. For that consider the  following  situation in $X_1$:
$$\xymatrix@=5pt{
	\alpha \ar[rr]^{R[d_1^{X_{\bullet}}]} \ar[dd]^{R[h_1]}  \ar@(l,l)[dd]_{R[d_0^{X_{\bullet}}]} && s_0d_1(\alpha) \ar[dd]_{R[h_1]} \ar@(r,r)@{.>}[dd]^{R[d_0^{X_{\bullet}}]}\\
	&&&\\
	\alpha' \ar[rr]_{R[d_1^{X_{\bullet}}]} && s_0d_1(\alpha')
}
$$
The left hand side bended arrow asserting that $d_0(\alpha)=d_0(\alpha')$ , the induced  right hand side dotted bended arrow means $d_1(\alpha)=d_1(\alpha')$. Now,  the pair $(h_1,d_1^{X_{\bullet}})$ being jointly monic by the pullback assumption, we get $\alpha=\alpha'$.
\endproof

\begin{coro}
	Let $\EE$ be a Gumm category and $(M,+,0)$ a (necessarily commutative) monoid in $\EE$.
	A morphism $\psi:M\times X \to X$ produces a left action of the monoid $M$ on $X$ as soon as $\psi(0,x)=x$, and this action is necessarily left cancellable.
\end{coro}

\section{Axiom $\ast$}\label{axiomst}

The aim now is to find a property which is fullfilled by any fiber $Pt_Y\VV$ when $\VV$ is a congruence modular variety, and which, mixed with the modular hyperextensible condition in the fibers $Pt_Y\EE$ of a Gumm category $\EE$, would allow us to achieve the proof of the construction of the associated abelian object in $Pt_Y\EE$ on the model of the construction given for the congruence modular varieties in Theorem \ref{main}. With this idea in mind, according to the proof in question, this property, or axiom, would be, in some way, related to factor permutation, but would preferably disymmetrize it in order to get a better understanding  of its behaviour.\\
\noindent\textbf{[Axiom $\ast$]:} \emph{Let $\EE$ be a pointed category; given any right punctual difunctional relation $W$ between $X$ and $Z$, and any equivalence relation $T$ on $W$:
$$\xymatrix@=10pt{
	W   \ar[dd]_{\pi_0^X}  \ar@<1ex>[rr]^{\pi_1^Z} &&  Z \ar[dd]\ar@{ >->}[ll]^{0W1}\\
	&&&&\\
	X  \ar@<1ex>[rr] && 1\ar@{ >->}[ll]^{0_X}  
}
$$
the forgetful canonical factorization $\zeta_{R[\pi_0^X],T}: R[\pi_0^X]\square T \to  R[\pi_0^X]\times_1T$ is an extremal epimorphism.} 

According to Proposition \ref{permutation} below, in a regular category $\EE$, the axiom $\ast$ means that, when the equivalence relations $R[\pi_0^X]$ and $R[\pi_1^Z]$ permute (difunctional relation), the equivalence relations $R[\pi_0^X]$ and $T$ do permute as well (extremal epimorphism). So, we get immediately:
\begin{prop}
	When a pointed regular category $\EE$ category  satisfies the Axiom $\ast$, it is factor permutable.
\end{prop}
\begin{coro}
	Let $\EE$ be any pointed regular category satisfying the Axiom $\ast$. Given any map $h: X\to Z$ and the cokernel $q$ of the monomorphism $(1_X,h): X\into X\times Z$, then the following induced downward square is regular pushout: 
	$$\xymatrix@=7pt{
		X\times Z    \ar[dd]_{p_0^X} \ar@{->>}[rrrr]^{q}  &&&& Q \ar[dd]_{} \\
		&&&&&&&\\
		X \ar@{ >->}@<-6pt>[uu]_{(1_X,h)} \ar@{->>}[rrrr]_{}  &&&& 1 \ar@<-6pt>@{ >->}[uu]_{0}
	}
	$$
	Taking $h=1_X$, we get the diagonal punctuation, and the following downward square is regular pushout:
	$$\xymatrix@=7pt{
		X\times X    \ar[dd]_{p_0^X} \ar@{->>}[rrrr]^{\omega_X}  &&&& Dp_X \ar[dd]_{} \\
		&&&&&&&\\
		X \ar@{ >->}@<-6pt>[uu]_{s_0^X} \ar@{->>}[rrrr]_{}  &&&& 1 \ar@<-6pt>@{ >->}[uu]_{0}
	}
	$$
\end{coro}
\proof
The relation $(p_0^X,p_1^Z)$ is right punctual and difunctional. Then considering the kernel equivalence relation $R[q]$ on $X\times Z$, the factorization $(p_0^X,q):X\times Z \to X\times Q$ is shown to be a regular epimorphism by the following permutation diagram:
$$\xymatrix@=5pt{
	(a,v) \ar@{.>}[rr]^{R[q]} \ar[dd]_{R[p_0^X]}   &&  (x,z) \ar@{.>}[dd]^{R[p_0^X]} \\
	&&&\\
	(a,h(a))\ar[rr]_{R[q]} && (x,h(x))
}
$$
\endproof
Our first goal is achieved with the following example:
\begin{prop}
	Let $\VV$ be any congruence modular variety. Then any fiber $Pt_Y\VV$ satisfies the Axion $\ast$.
\end{prop}
\proof
Consider any right punctual difunctional relation:
$$\xymatrix@=10pt{
	W   \ar[dd]_{\pi_0^X}  \ar@<1ex>[rr]^{\pi_1^Z} &&  Z \ar[dd]_{\phi_Z}\ar@{ >->}[ll]^{\sigma \phi W1_Z}\\
	&&&&\\
	X  \ar@<1ex>[rr]^{\phi_X} && Y\ar@{ >->}[ll]^{\sigma_X}  \ar@<-1ex>@{ >.>}[uu]_{\sigma_Z}
}
$$
From the existence of the map  $\sigma \phi W1_Z$, the fact that the relation $W$ is difunctional insures that, if we have $aWb$, we then get $aWc$ for any $c\in \phi_Z^{-1}(\phi_X(a))$. 
Let $T$ be an equivalence relation on $W$ in the fiber $Pt_Y\VV$. This means that $T$ is annihilated by $\phi_X.\pi_0^X$. Now suppose that $(aWc)T(xWz)$; thus we get $\phi_Z(c)=\phi_X(a)=\phi_X(x)=\phi_Z(z)$. Then, from $aWb$, consider the following situation depicted by the plain arrows:
$$\xymatrix@=10pt{
	&&	aWb \ar@{.>}[lldd]_{T}  \ar[dd]_{\pi_1^Z}  \ar[rr]^{\pi_0^X} &&  aWz \ar[dd]_{\pi_1^Z} \ar[rr]^{\pi_0^X} && aWc \ar[dd]^{\pi_1^Z} \ar[lldd]_T\\
	&&&&&&\\
	t(xWb,xWc,xWz)\ar@{.>}[rr]_>>>>>>{\pi_0^X}  && xWb  \ar[rr]_{\pi_0^X} && xWz \ar[rr]_{\pi_0^X} && xWc 
}
$$
Through the ternary term $t$ on $W$ given  by Theorem 4.3 in \cite{Gu}, from these plain arrows, we get the dotted ones which insure the permutation of the equivalence relations $R[\pi_0^X]$ and $T$. 
\endproof
And our main goal is achieved with the following:
\begin{theo}\label{main2}
Let $\EE$ be a regular Gumm category such that any fiber $Pt_Y\EE$ satisfies the Axion $\ast$. If, in addition, these fibers admit diagonal punctuations, then, given any split epimorphism $(f,s): X\splito Y$, the diagonal punctuation in $Pt_Y\EE$:
	$$\xymatrix@=7pt{
	R[f]  \ar@<-6pt>@{.>}[dd]_{d_0^f}  \ar@{->>}[rrr]^{\omega_f}  &&& Dp[f] \ar@<6pt>@{.>}[dd]^{\psi_f} \\
	&&&&&&\\
	X \ar@{ >->}[uu]_{s_0^f} \ar@{->>}[rrr]_{f}  &&& Y \ar@{ >->}[uu]^{\theta_f}
}
$$	
produces with $(\psi_f,\theta_f)$ the universal abelian split epimorphism associated with $(f,s)$. The subcategory $AbPt_Y\EE\subset Pt_Y\EE$ is stable under monomorphism. 
\end{theo}
\proof
We can now follow step by step the proof of Theorem \ref{main}. Consider the following diagram:
$$\xymatrix@=10pt{
	R[\omega_{f}] \ar@<-6pt>[rrr]_{d_0^{\omega}} \ar@<6pt>[rrr]^{d_1^{\omega}}
	\ar@<-6pt>[dd]_{R(d_0^f)}   \ar@<6pt>[dd]^{R(d_1^f)}	&&& R[f] \ar@{ >->}[lll] \ar@<-6pt>[dd]_{d_0^f}   \ar@<6pt>[dd]^{d_1^f} \ar@{->>}[rrrr]^{\omega_f}  &&&& Dp[f] \ar@<6pt>@{.>}[dd]^{\psi_f} \\
	&&&&&&&\\
	R[f] \ar@<-6pt>[rrr]_{d_0^f}   \ar@<6pt>[rrr]^{d_1^f} \ar@{ >->}[uu]	&&& X \ar@{ >->}[uu]|{s_0^f} \ar@{->>}[rrrr]_{f} \ar@{ >->}[lll] &&&& Y \ar@{ >->}[uu]^{\theta_f}
}
$$
The category $\EE$ being a Gumm category, the left hand side vertical reflexive relation is an equivalence relation by Theorem \ref{DD}; and the right hand side square is regular pushout by the previous proposition since the fiber $Pt_Y\EE$ satisfies the Axiom $\ast$.
\endproof
Here is now another easy varietal example of pointed categories satisfying the Axiom $\ast$.
\begin{prop}
Let $\VV$ be a strongly unital variety. Then the variety $\VV$ satisfies the Axiom $\ast$.
\end{prop}
\proof
Let $p$ be the ternary term insuring the strong unitality with $p(x,x,y)=y$ and $p(x,0,0)=x$. Consider any  right punctual difunctional relation:
 $$\xymatrix@=10pt{
 	W   \ar[dd]_{\pi_0^X}  \ar@<1ex>[rr]^{\pi_1^Z} &&  Z \ar[dd]\ar@{ >->}[ll]^{0W1}\\
 	&&&&\\
 	X  \ar@<1ex>[rr] && 1\ar@{ >->}[ll]^{0_X}  
 }
 $$
 and any equivalence relation $T$ on $W$.  Starting with the plain arrows, we get: 
 $$\xymatrix@=5pt{
 	aWb \ar@{.>}[rr]^{T} \ar[dd]_{R[\pi_0^X]}   &&  xWp(y,c,b) \ar@{.>}[dd]^{R[\pi_0^X]} \\
 	&&&\\
 	aWc\ar[rr]_{T} && xWy
 }
 $$
by $aWb=(aWc,0Wc,0Wb)Tp(xWy,0Wc,0Wb)=(x,p(y,c,b))$. 
\endproof

The weaker property of pure factor permutability for such a kind of varieties was already noticed in \cite{Gfp}.

\section{Congruence permutation and connector}

This section is mainly devoted to recalls and precisions about the relationship between permutation and centralization of pairs $(R,S)$ of equivalence relations on an object $X$ in any category  $\EE$.

\subsection{The square construction}\label{square}

Given a pair $(R,S)$ of equivalence relations on an object $X$ in a  category $\mathbb E$, we denote by $R\square S$ the inverse image of the equivalence relation $S\times S$ on $X\times X$ along the inclusion $(d_0^R,d_1^R):R\rightarrowtail X\times X$. This defines a double equivalence relation:
$$ \xymatrix@C=3pc@R=2pc{
	R \square S \ar@<-1,ex>[d]_{\delta_0^R}\ar@<+1,ex>[d]^{\delta_1^R} \ar@<-1,ex>[r]_{\delta_0^S}\ar@<+1,ex>[r]^{\delta_1^S}
	& S \ar@<-1,ex>[d]_{d_0^S}\ar@<+1,ex>[d]^{d_1^S} \ar[l]\\
	R \ar@<-1,ex>[r]_{d_0^R} \ar@<+1,ex>[r]^{d_1^R} \ar[u]_{} & X\ar[u]_{} \ar[l]
}
$$
In set theoretical terms, $R\square S$ is the set of quadruples $(x,y,t,z)\in X^4$ such that $xRy$, $tRz$, $xSt$ and $ySz$, often depicted as:
$$\xymatrix@C=1,5pc@R=1,5pc{x \ar@{-}[r]^-{S} \ar@{-}[d]_-{R} & t \ar@{-}[d]^-{R} \\ y \ar@{-}[r]_{S} & z}$$ 
A very useful tool is given by the following forgetful factorization:
\begin{equation} \label{canonical}
\zeta_{R,S}\colon R \square S \to R \times_1 S\colon (x,y,t,z) \mapsto xRySz.
\end{equation}
where $R \times_1 S$ is defined by the following pullback:
$$\xymatrix@=4pt{
&&&&&&&&&&& R \times_1 S  \ar@<-1ex>[ddd]_{\bar d_0^R}   \ar@<1ex>[rrrrr]^{\bar d_1^S}  &&&&&  S  \ar@<-1ex>[ddd]_{d_0^S} \ar[lllll]^{\bar s_0^S}
&&&&& && &&&&&&&&\\
&&& &&&&&&&&&&&&&\\
&&& &&&&&&&&&&&&&\\
&&&&&&&&&&& R \ar[uuu]_{\bar s_0^R}  \ar@<1ex>[rrrrr]^{d_1^R}  &&&&& X \ar[lllll]^{s_0^R} \ar[uuu]_{s_0^S}
}
$$
If $R \cap S = \Delta_X$ (the diagonal equivalence relation on $X$), this factorization $\zeta_{R,S}$ is a monomorphism. Indeed, if $(x,y,t,z)$ and $(x,y,t',z)$ are in $R\square S$, then $tRzRt'$ and $tSxSt'$, showing that $t (R\cap S) t'$ and thus $t=t'$.
\begin{defi}
We shall say that a pair $(R,S)$ of equivalence relations on $X$ has a \emph{sharp} intersection when $R\cap S=\Delta_X$ and when, in addition, the monomorphism $\zeta_{R,S}$ is an isomorphism. A relation $(f,g): W \into X\times Z$ is \emph{difunctional} when the pair $(R[f],R[g])$ has a sharp intersection. 
\end{defi}
When $R$ and $S$ have a sharp intersection, then any commutative square determined by the diagram of the double equivalence relation associated with $R\square S$ is a pullback. 
On the other hand, it is well known that a reflexive relation is an equivalence relation if and only if it is difunctional, see  \cite{CLP} for instance. Straighforward is the following observation:
\begin{lemma}
	Given any pair $(R,S)$ with a sharp intersection in a category $\EE$, the following diagram is underlying an equivalence relation:
	$$ \xymatrix@C=3pc@R=2pc{
		R \square S  \ar@<-1,ex>[r]_{d_1^R.\delta_0^R}\ar@<+1,ex>[r]^{d_1^S.\delta_0^S}
	  & X \ar[l]
	}
	$$
	which is the supremum $R\vee S$ of $R$ and $S$ among the equivalence relations on $X$. Moreover we get: $R\circ S=R\vee S=S\circ R$ without any regular assumption on the category $\EE$.
\end{lemma}
Straightforward as well is the following:
\begin{lemma}\label{keyperm}
	Given any pair $(R,S)$ of equivalence relations  on $X$ in a category  $\EE$, the following commutative square is a pullback in $\EE$:
	$$ \xymatrix@C=3pc@R=2pc{
		R \square S \ar[d]_{\zeta_{S,R}} \ar[r]^{\zeta_{R,S}}
		& R\times_1S \ar[d]^{(d_0^R.\bar d_0^R,d_1^S.\bar d_1^S)}\\
		S\times_1R \ar[r]_{(d_0^S.\bar d_0^S,d_1^R.\bar d_1^R)}   & X\times X
	}
	$$
\end{lemma}

\subsection{Congruence permutation}\label{conn}

In this section we shall suppose $\EE$ regular, so that, in any case, we can compose the relations in $\EE$. The composition $S\circ R$ of two equivalence relations on $X$ is given by the canonical decomposition of the above vertical right hand side map: $$(d_0^R.\bar d_0^R,d_1^S.\bar d_1^S): R\times_1S \onto S\circ R \into X\times X$$  We have $S\circ R\subset R\circ S$, when we get the following dotted factorization making commutative the lower right hand side triangle and, consequently, the vertical rectangle as well:
	$$ \xymatrix@C=3pc@R=2pc{
	R \square S \ar[dd]_{\zeta_{S,R}} \ar[r]^{\zeta_{R,S}}
	& R\times_1S \ar@{->>}[d]_{\bar{\rho}} \ar[ddr]^{(d_0^R.\bar d_0^R,d_1^S.\bar d_1^S)}\\
	& S\circ R \ar@{ >.>}[d] \ar@{ >->}[dr]^{\bar j} \\
	S\times_1R \ar@<-2ex>[rr]_{(d_0^S.\bar d_0^S,d_1^R.\bar d_1^R)}\ar@{->>}[r]^{\rho}   & R\circ S \ar@{ >->}[r]^j & X\times X
}
$$

\begin{prop}\label{permutation}
Let $\EE$ be any regular category. Given any pair $(R,S)$ of equivalence relations on an object $X$, the following conditions  are equivalent:\\
1) $S\circ R = R\circ S$;\\
2) $S\circ R \subset R\circ S$;\\
3) the above vertical rectangle is a pullback;\\
4) $\zeta_{R,S}$ is a regular epimorphism;\\
5) $\zeta_{S,R}$ is a regular epimorphism.
\end{prop}
\proof
When $R$ and $S$ are equivalence relations, then it is clear that $S\circ R \subset R\circ S$ if and  only if $R\circ S\subset S\circ R$, namely $(1)\iff (2)$, since $(S\circ R)^{op}=R^{op}\circ S^{op}$.

Suppose $(2)$. The above whole quadrangle being  a pullback by the previous lemma, and the map $j$ being a monomorphism the vertical rectangle is a pullback; whence $(3)$. Suppose $(3)$, the map $\rho$ being a regular epimorphism, so is $\zeta_{R,S}$; whence $(4)$. Suppose $(4)$. Then the map $\bar{\rho}.\zeta_{R,S}$ is a regular epimorphism. So, the presence of the monomorphism $j$ in the commutative quadrangle induces the dotted factorization, since, in any regular category, regular epimorphisms and strong epimorphisms coincide. Whence $(2)$. We get $(1) \iff (5)$ by inverting the role of $R$ and $S$.
\endproof

A good way to get a regular epimorphism is to produce a splitting, whence the following:
\begin{defi}
Given a category $\EE$ and a pair $(R,S)$ of equivalence relations on an object $X$, a \emph{preconnector} for this pair is a splitting of $\zeta_{R,S}$. 
\end{defi}

\subsection{Recalls about connectors}

Connectors between pairs $(R,S)$ of equivalence relations were introduced in \cite{BG} mainly to investigate the categorical aspect of the \emph{centralization} of equivalence realtions  in the Mal'tsev context:
\begin{defi}\label{connector}
Let $R$ and $S$ be two equivalence relations on a same object $X$ in $\EE$, a \emph{connector} between $R$ and $S$ is a morphism $p\colon R \times_X S \to X$ satisfying the following axioms:
\begin{enumerate}
	\item $xSp(xRySz)$ and $p(xRySz)Rz$;
		$$\xymatrix@!0@R=3em@C=3em{x \ar@{-}[rr]^{R} \ar@{-}[rd]_-{S} && y \ar@{-}[rd]^-{S} \\ & p(xRySz) \ar@{-}[rr]_-{R} && z}$$
	\item $p$ is a partial Mal'tsev operation, i.e., $p(xRxSy)=y$ and $p(xRySy)=x$;
	\item $p$ is left and right associative, i.e., $p(p(xRySz)RzSw)=p(xRySw)$ and $p(xRySp(yRzSw))=p(xRzSw)$.
\end{enumerate}
\end{defi}
So, in a regular context, a connector $p$ is nothing but a way to produce a \emph{coherent permutation} $R\circ S=S\circ R$ which is denoted by $[R,S]_p=0$. Clearly:\\
1) Axiom $1)$ produces a splitting of $\zeta_{R,S}$, and thus a preconnector we shall denote by $\pi_p$;\\
2) a pair $(R,S)$ of equivalence relation with a sharp intersection is connected, in a unique way, via the inverse of the isomorphic comparison $\zeta_{R,S}: R\square S \to R\times_1S$;\\
3) the notion of pseudogroupoid introduced in \cite{JP} shows that given any reflexive graph:
$$\xymatrix@=22pt
{
	G \ar@<1ex>[rr]^{d_{1}} \ar@<-1ex>[rr]_{d_{0}} &&
	X_0   \ar[ll]|{s_0}  
}
$$
there is bijection between the groupoid structures on this reflexive graph and the connectors between the kernel equivalence relations $R[d_0]$ and $R[d_1]$.

In \cite{BG}, it was shown that in a Mal'tsev category, namely a category in which any reflexive relation is an equivalence relation \cite{CLP} \cite{CPP}, Axiom $2)$ implies the two others. In a Gumm category, we observed here, and previously in \cite{BGGGum}, that we need Axioms $1)$ and $2)$ to imply the third one.

\medskip

\noindent {\bf keywords}: Mal'tsev and congruence modular variety; Mal'tsev and Gumm category; split epimorphism; internal category and groupoid; J\'onnson-Tarski and subtractive variety; unital, subtractive  and strongly unital categories.\\
{\bf Mathematics Subject Classification (2020).} 08A30, 08B05, 18A20, 18C10, 18C40, 18E13.

\medskip

\noindent Univ. Littoral C\^ote d'Opale, UR 2597, LMPA,\\
Laboratoire de Math\'ematiques Pures et Appliqu\'ees Joseph Liouville,\\
F-62100 Calais, France.\\
bourn@univ-littoral.fr


\begin{thebibliography}{99}
	
	\bibitem{Barr} M.\ Barr, \emph{Exact categories}, Springer L.N. in Math. \textbf{236} (1971), 1--120.
	
	\bibitem{BB} F.\ Borceux and D.\ Bourn, Mal'cev, Protomodular, Homological and Semi-Abelian Categories, \textit{Kluwer, Mathematics and Its Applications} \textbf{566} (2004), 479 pp.
	
	\bibitem{B0} D.\ Bourn, \emph{Normalization equivalence, kernel equivalence and affine categories}, Springer L.N. in Math. \textbf{1488} (1991), 43--62.
	
	\bibitem{Bfib} D.\ Bourn, \emph{Mal'cev Categories and fibration of pointed objects}, Applied categorical structures, \textbf{4} (1996), 302--327.
	
	\bibitem{B3} D.\ Bourn, \emph{Baer sums and fibered aspects of Mal'cev operations}, Cahiers de Top. et G\'eom. Diff. \textbf{40} (1999), 297-316.
	
	\bibitem{Bint} D.\ Bourn, \emph{Intrinsic centrality and associated classifying properties}, Journal of Algebra, \textbf{256} (2002), 126-145.
	
	\bibitem{BGum} D.\ Bourn, \emph{Fibration of points and congruence modularity}, Algebra Universalis \textbf{52} (2005), 403-429.
	
	\bibitem{Bab} D.\ Bourn, \emph{Abelian groupoids and non-pointed additive categories},  Theory and Applications of Categories, \textbf{20} (2008), 48-73.
	
	\bibitem{Bmon} D.\ Bourn, \emph{On the monad of internal groupoids}, Theory and Applications of Categories, \textbf{28} (2013), 150-165.
	
	\bibitem{BG} D.\ Bourn and M.\ Gran, \emph{Centrality and connectors in Maltsev categories}, Algebra univers. \textbf{48} (2002), 317-343.
	
	\bibitem{BGGGum} D.\ Bourn and M.\ Gran, \emph{Categorical aspects of modularity}, in :Galois theory, Hopf algebras, and Semiabelian categories, G.Janelidze, B.Pareigis, W.Tholen editors,  Fields Institute Communications, \textbf{43}, Amer. Math. Soc. (2004), 77-100.
	
	\bibitem{BGGum} D.\ Bourn and M.\ Gran, \emph{Normal sections and direct product decompositions}, Communications in Algebra, \textbf{32} (2004), 3825-3842.
	
	\bibitem{BZZ} D.\ Bourn and Z.\ Janelidze, \emph{Subtractive categories and extended subtractions}, Applied categorical structures \textbf{17} (2009), 302-327.
	
	\bibitem{BZ} D.\ Bourn and Z.\ Janelidze, \emph{A note on the abelianization functor}, Communications in Algebra \textbf{44} (2016), 2009-2033.
	
	\bibitem{CLP} A.\ Carboni, J.\ Lambek and M.C.\ Pedicchio, \emph{Diagram chasing in Mal'cev categories}, J. Pure Appl. Algebra \textbf{69} (1991), 271-284.
	
	\bibitem{CPP} A.\ Carboni, M.C.\ Pedicchio and N.\ Pirovano, \emph{Internal graphs and internal groupoids in Mal'cev categories}, CMS Conference Proceedings \textbf{13} (1992), 97-109. 
	
	\bibitem{Gfp} M.\ Gran, \emph{Applications of Categorical Galois Theory in Universal Algebra}, in :Galois theory, Hopf algebras, and Semiabelian categories, G.Janelidze, B.Pareigis, W.Tholen editors,  Fields Institute Communications, \textbf{43}, Amer. Math. Soc. (2004), 243-280.
	
	\bibitem{Gu} H.P.\ Gumm, \emph{Geometrical methods in congruence modular varieties}, Mem. Amer. Math. Soc. \textbf{45} (1983).
	
	\bibitem{JP} G.\ Janelidze and M.C.\ Pedicchio, \emph{Pseudogroupoids and commutators}, Theory and Applications of Categories \textbf{8} (2001), 408-456.
	
	\bibitem{ZJ} Z.\  Janelidze, \emph{Subtractive categories}, Applied categorical structures \textbf{13} (2005), 343-350.
	
	\bibitem{J} P.T.\ Johnstone, Topos Theory, Academic Press, London, (1977).
	
	\bibitem{Jo} P.T.\ Johnstone, \emph{Affine categories and naturally Mal'cev categories}, Journal of Pure and Applied Algebra, \textbf{61} (1989), 251--256.
	
	\bibitem{Ta} A.\ Tarski, \emph{Ein Beitrag zur Axiomatik der Abelschen Gruppen}, Fund. Math. \textbf{30} (1938), 253--256.
	
	\bibitem{Ur} A.\ Ursini, \emph{0n subtractive varieties}, Algebra univers. \textbf{31} (1994), 204-222.
\end{thebibliography}
\end{document}